\documentclass[11pt,twoside]{article}

\usepackage{epsfig,amsfonts,color}
\usepackage{amsmath}
\usepackage{mathtools}

\usepackage{amssymb, palatino, geometry,url}
\usepackage{mathrsfs}
\usepackage{amsthm}
\usepackage{algorithm2e}

\usepackage{booktabs}
\usepackage{siunitx}
\usepackage{xcolor}
\usepackage{url}
\usepackage[most]{tcolorbox}
\usepackage[T1]{fontenc}
\newtcolorbox{highlighted}{colback=cyan,coltext=black,breakable}

\usepackage[colorlinks=true,linkcolor=black,citecolor=black,urlcolor=black]{hyperref}
\usepackage{subcaption}
\usepackage{cleveref}
\usepackage{color}
\usepackage{fancyhdr}
\theoremstyle{definition}

\newcommand{\be}{\begin{equation}}
\newcommand{\ee}{\end{equation}}
\newcommand{\bea}{\begin{eqnarray}}
\newcommand{\eea}{\end{eqnarray}}

\newcommand{\bvec}{\left(\begin{array}{c}}
\newcommand{\evec}{\end{array}\right)}
\newcommand{\bsub}{\begin{subequations}}
\newcommand{\esub}{\end{subequations}}

\usepackage[thinc]{esdiff}

\usepackage{lineno}
\usepackage{tabularx}
\usepackage{multirow}
\usepackage{adjustbox}
\usepackage[outdir=./Output]{epstopdf}

\begin{document}

\title{Degradation-Aware Model Predictive Control for Battery Swapping Stations under Energy Arbitrage}

\author{
Ruochen Li\textsuperscript{1}, Zhichao Chen\textsuperscript{1}, Zhaoting Zhang \textsuperscript{1},
Renjie Guo\textsuperscript{2}, Zhankun Sun\textsuperscript{2}, Jiwei Yao\textsuperscript{3}, Jiaze Ma\textsuperscript{1}\thanks{Corresponding author: jiazema@cityu.edu.hk} \\
\\
\small\textsuperscript{1}Department of Systems Engineering, City University of Hong Kong, Kowloon, Hong Kong, China \\
\small\textsuperscript{2}Department of Decision Analytics and Operations, City University of Hong Kong, Kowloon, Hong Kong, China \\
\small\textsuperscript{3}Department of Chemical Engineering, University of Utah, Salt Lake City, UT, USA
}
\date{}
\maketitle

\begin{abstract}
Battery swapping stations (BSS) offer a fast and scalable alternative to conventional electric vehicle (EV) charging, gaining growing policy support worldwide. However, existing BSS control strategies typically rely on heuristics or low-fidelity degradation models, limiting profitability and service level. This paper proposes BSS-MPC: a real-time, degradation-aware Model Predictive Control (MPC) framework for BSS operations to trade off economic incentives from energy market arbitrage and long-term battery degradation effects. BSS-MPC integrates a high-fidelity, physics-informed battery aging model that accurately predicts the degradation level and the remaining capacity of battery packs. The resulting multiscale optimization—jointly considering energy arbitrage, swapping logistics, and battery health—is formulated as a mixed-integer optimal control problem and solved with tailored algorithms. Simulation results show that BSS-MPC outperforms rule-based and low-fidelity baselines, achieving lower energy cost, reduced capacity fade, and strict satisfaction of EV swapping demands. 
\end{abstract}


\section{Introduction}

The global shift toward electric vehicle (EVs) is accelerating. Yet, despite progress in battery and charging technologies, two obstacles continue to limit EV adoption: range anxiety and charging bottlenecks. Drivers remain concerned about running out of power during long trips, especially on highways. While in densely populated urban areas where people lack private charging space, peak charging demand (e.g. during holidays) often results in long queues and hours of waiting. For high-utilization EV such as taxis, buses, and trucks, every minute of downtime directly translates into lost productivity~\cite{Millard2023BatterySwapping}.

BSS offer a solution to these problems by allowing EV drivers to replace a depleted battery with a fully charged one in just a few minutes. A well-known example is NIO’s BSS, which consists of two integrated components: a garage that robotically swaps the battery packs and a charging stack that holds and recharges 21 depleted batteries \cite{nio2025apparticle527488}. NIO’s fully automated stations can complete a battery exchange in about three minutes \cite{nio_power}. Compared with fast charging, BSS eliminates long wait times while avoiding the heavy grid loads and accelerated battery degradation caused by fast charging. Beyond convenience, BSS also decouple battery ownership from vehicle ownership, which can reduce the upfront cost of EVs for consumers. Centralized battery handling also enables efficient recycling and lifecycle management of EV battery \cite{DOE2022BatteryReport}. 

Recent years have witnessed substantial policy and industrial momentum for battery swapping. For example, NIO has deployed more than 3500 BSS cross China and provided battery swapping service over 87 million times~\cite{NIO2024SwapStations}. Ecosystem leaders like Gogoro~\cite{Gogoro2024SwapNetwork} and Ample~\cite{FT2025BatterySwappingVC} have initiated swap networks for both two- and four-wheelers. Major policy initiatives in China~\cite{MIIT2020NEVPlan} and pilot programs in the US and Europe~\cite{YangTang2024BatterySwapForbes} reflect a growing consensus that integrating battery swapping with smart grid management accelerates sustainable mobility.

Despite demonstrated technical feasibility—most visibly in Tesla’s early battery-swap trials—market adoption has lagged because economics, not engineering, are the binding constraint~\cite{Loveday2015ElonMuskBatterySwap,zhan2022review}. A BSS requires upfront investment and ongoing maintenance, while consumers’ willingness to pay service fees is limited~\cite{wu2021survey}. Pricing high erodes demand; pricing low impedes cost recovery—producing a profitability dilemma that has hindered large-scale deployment.

A promising path forward is to treat the BSS not only as a mobility service but also as a grid-integrated energy storage asset. By adopting operational modes that co-optimize swapping with energy arbitrage—charging during low-price periods and discharging (or offsetting purchases) during high-price periods—the BSS can augment revenue and/or reduce electricity cost. This dual role reframes station operation from cost center to value-creating flexibility resource. For example, NIO’s newly deployed BSS can store up to 21 EV batteries with a total capacity of 2.1 MWh. Leveraging this energy reservoir creates an additional revenue stream that enhances the economic viability of the stations. Specifically, advanced operation can unlock three key benefits\cite{cui2023operation,revankar2021grid,wu2021survey}:

\begin{itemize}
  \item \textbf{Demand-driven preparation:} Given the swapping demand over a short horizon, the station can proactively charge selected batteries in advance. It ensures sufficient availability of fully charged batteries at all times while avoiding last-minute rushes.
  
  \item \textbf{Market-aware energy scheduling:} Given the electricity price trends, the system can strategically shift charging to low-price periods and discharge during high-price intervals, thereby maximizing profit and enhancing participation in energy markets.
  
  \item \textbf{Degradation-aware operation:} Since battery aging is influenced by usage patterns, state-of-charge levels, and current profiles, integrating battery degradation models into control decisions allows the station to extend battery life. This mitigates long-term economic loss casued by the degradation of EV battery.
\end{itemize}

The above factors position BSS as a systems that must balance real-time service reliability, long-term asset health, and economics. Yet implementing predictive, health-aware control in real time is challenging: high-fidelity degradation models entail nonlinear dynamics and high-dimensional state spaces, leading to  computational burden \cite{cao2020multiscale}. Conversely, simplified models and heuristic battery-management rules constrain the station’s ability to fully exploit its value as a grid-integrated energy storage resource.

To overcome these challenges, this paper proposes a MPC framework of BSS with integrated energy market arbitrage. By incorporating a high-fidelity battery degradation model via a computationally efficient surrogate, the framework enables real-time, degradation-aware decision-making across both battery swapping logistics and grid-level energy transactions. This model-based strategy not only ensures short-term operational efficiency but also promotes long-term battery health.

\noindent\textbf{Organization:} The remainder of this manuscript is organized as follows. ~\Cref{sec:RelatedWorks} provides a review of existing studies and identifies the open challenges that still need to be resolved. Building on these insights, our proposed control framework is detailed in~\Cref{sec:proposedApproachSection}. Finally, in Section~\Cref{sec:Experiment}, we first present the results of overall charging and discharging within a prediction horizon, illustrating the response of BSS-MPC to price and demand fluctuations. This part includes a comparison between the planning results at the initial time step and the results obtained from rolling horizon optimization, showing how the control inputs are adjusted and gradually differ as the sliding window progresses in ~\cref{sec:planning_result}. In the second part of the experiment ~\cref{sec:comparison}, we assess the effectiveness of BSS-MPC by comparing it with other methods, including rule-based and simplified model baselines.

\section{Literature Review}\label{sec:RelatedWorks}

Early studies used either rule-based or simplified battery model to schedule charging, discharging and swapping in BSS. The models in \cite{YANG2014544,mahoor2017electric, 8054726} simply classified battery states into fully charged, half-charged, half-discharged, or empty to enable a straightforward management strategy. The rule-based BSS scheduling model in \cite{6590033} charges only during off-peak hours to minimize cost, permits but deprioritizes mid-peak charging, and avoids peak-hour charging entirely. The model in \cite{7165769} uses empirical formulas to approximate battery charging status without considering the detailed battery electrochemistry, relying instead on a simplified constant current-voltage estimation model. To account for battery evolution characteristics and exploit grid arbitrage, subsequent studies on BSS operation moved beyond heuristic and rule-based strategies. In \cite{wu2017charging}, the model in \cite{7165769} is extended with a surrogate for battery aging. The objective is to maintain a high stock of swappable batteries (service reliability) while minimizing charging-induced degradation; However, long-term arbitrage opportunities are not considered. \cite{sarker2013electric,7285624} represent state-of-charge (SOC) dynamics as a linear function of charging and discharging power, while battery degradation was approximated—also linearly—as a function of charging and discharging power. The problem is cast as an optimal control and solved in a rolling-horizon manner. However, these simplifications —such as discretizing charging/discharging modes and linearizing nonlinear effects-induce cumulative errors and degrade forecasting accuracy over long horizons. 

Although MPC gradually emerged as the mainstream framework for BSS scheduling \cite{infante2019optimal,10506804,ahmad2019cost}, most implementations deliver implementable rules but typically on simplified battery state models, again curbing fidelity for aging and long-run arbitrage valuation. Parallel works \cite{wu2021two,esmaeili2019optimal,vsepetanc2019cluster} from the market perspective frames the BSS as a price-taker optimizing participation across day ahead /real time/ ancillary services with demand forecasts are useful to expose multi-market value streams, but most studies either omit explicit electrochemical aging or approximate it coarsely, which biases long-horizon economics. 

To address computational challenges from complex electrochemistry and combinatorics, recent learning-based studies use reinforcement learning (RL) for fast online decisions, sometimes paired with a post-hoc optimization layer—for example, a DRL–optimization cascade for scalable BSS energy management \cite{su2024energy} and RL policies with degradation/waiting-time penalties in the reward \cite{shalaby2023model}. While these approaches scale, they typically encode hard constraints and battery health only indirectly (via penalties or post-processing), or rely heavily on offline training that makes performance sensitive to model assumptions, and they provide limited optimality guarantees. Consequently, they still fall short of solving the full mixed-integer, health-aware arbitrage problem at high fidelity. Decomposition schemes likewise highlight the computational burden and depend on model simplifications to remain tractable \cite{sohrabi2024electrification}.

We observed that two gaps recur across the literature:
\begin{enumerate}
\item Battery-model fidelity: Aging physics and its feedback on efficiency and usable capacity are often simplified, degrading the valuation of long-term arbitrage and risking inventory shortfalls; 
\item Tractability: Health-aware, high-fidelity scheduling couples nonlinear battery dynamics with integer swap/charge decisions, making the problem computationally hard; most works therefore revert to simplified states.
\end{enumerate}

There is a growing body of work showing that high-fidelity, degradation-aware models can materially improve stationary battery arbitrage decisions. Reniers et al. \cite{RENIERS201891} compare physics-based model and low-fidelity surrogates over year-long arbitrage and show substantially better revenue–lifetime valuation with the high fidelity physics-based model. Building on this, Cao et al. \cite{cao2020multiscale} embed a physics-based battery in a multiscale MPC with a terminal fade penalty, demonstrating closed-loop feasibility and large lifetime/profit gains versus low-fidelity MPC. Complementing these, Perez et al. \cite{7525538} solve optimal control directly on an electrochemical-thermal model with state constraints, illustrating how high-order dynamics can be handled numerically and how constraints map to practical charging/dispatch policies. These works indicate that high-fidelity control is tractable and economically meaningful for stationary storage—motivating our BSS extension, where integer swap logistics and service-level constraints introduce additional combinatorial complexity.

To address the aforementioned challenges, this paper presents, to the best of our knowledge, the first computational framework that systematically integrates a high-precision, physics-informed, degradation-aware model with MPC for real-time operation of BSS participating in energy market arbitrage. The main contributions are summarized as follows:
\begin{enumerate}
\item We formulate the BSS operation problem as a finite-horizon optimal control problem within a multiscale framework. This formulation explicitly captures the coupling between short-term energy dispatch decisions and long-term battery degradation dynamics under volatile electricity market signals.
\item To address the computational difficulties brought by the high-fidelity battery degradation model, we replace the original Single Particle Model (SPM) with a Kriging-based surrogate model that preserves the high-fidelity characteristics of both battery degradation behavior and charge/discharge dynamics, while largely enhancing computational efficiency and numerical stability. This enables the integration of battery aging effects into real-time MPC.
\item We solve the resulting optimal control problem as a mixed-integer nonlinear programming (MINLP) problem and develop a tailored MPC framework, namely \emph{BSS-MPC}. This framework leverages updated surrogate models and receding-horizon optimization to achieve coordinated decision-making with improved economic performance and enhanced battery longevity.
\end{enumerate}

\section{Proposed Approach}\label{sec:proposedApproachSection}
\begin{figure*}[!h]
    \centering
    \includegraphics[width=1\linewidth]{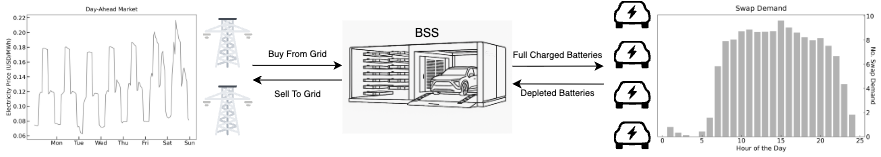}
    \caption{Schematic of BSS management with energy market arbitrage.}
    \label{fig:overall_workflow_result}
\end{figure*}

The application context of the proposed BSS-MPC framework is depicted in ~\Cref{fig:overall_workflow_result}.  We assume that the system operates within a Day-Ahead Market of PJM \cite{PJMMarketsOperations}, where electricity prices and battery swapping demands are known in advance and vary on an hourly basis (we consider a deterministic process). Based on energy prices and demands over a future time horizon, the BSS makes optimal scheduling decisions using high-fidelity degradation and charging models. Additionally, we assume that batteries supplied to EVs always operate within a specified SOC working range under stable discharge conditions, such that battery life degradation due to EV usage can be neglected. The goal is to maximize operational profit by strategically purchasing or selling electricity to charge batteries for future swapping, thereby enabling energy arbitrage. Meanwhile, the system must ensure that sufficiently charged batteries are always available to meet real-time swapping requests from EV drivers.

This chapter is organized as follows:
\begin{itemize}
    \item \Cref{sec:SPM} introduces the high-fidelity electrochemical battery model, including both charge/discharge dynamics and degradation behavior.
    \item \Cref{sec:Surrogate} develops a Kriging-based surrogate model that approximates the high-fidelity model while alleviating computational burdens of the optimization problem.
    \item \Cref{sec:Multiscale_Optimization_of_BSS} formulates the operation of BSS as a finite-horizon mixed-integer optimal control problem.
    \item \Cref{sec:BSS-MPC} presents the proposed BSS-MPC framework, which solves the resulting MINLP problem in a receding-horizon fashion.
\end{itemize}

%


\subsection{Single Particle Model} \label{sec:SPM}
\begin{table*}[!h]
\centering
\caption{Variables and Parameters in the Single Particle Model}
\label{tab:sp_variables}
\begin{tabularx}{\textwidth}{lXll}
\toprule
\textbf{Symbol} & \textbf{Description} & \textbf{Unit} & \textbf{Type} \\
\midrule
$C_{p}^{\text{avg}}$, $C_{n}^{\text{avg}}$ & Average lithium concentration (positive/negative electrodes) & $\mathrm{mol}/\mathrm{m}^3$ & State \\
$\delta_{\text{SEI}}$ & SEI (Solid Electrolyte Interphase) layer thickness & $\mathrm{m}$ & State \\
$c_f$ & Cumulative capacity fade & $\mathrm{Ah}$ & State \\
$P$ & Applied power (positive: discharging to grid; negative: charging from grid) & $\mathrm{MW}$ & Input \\
$C_p^{\text{surf}}$, $C_n^{\text{surf}}$ & Surface lithium concentration (positive/negative electrodes) & $\mathrm{mol}/\mathrm{m}^3$ & Algebraic \\
$i_{\text{int}}$ & Internal current & $\mathrm{A}$ & Algebraic \\
$\phi_p$, $\phi_n$ & Solid-phase potentials in positive/negative electrodes & $\mathrm{V}$ & Algebraic \\
$pot$ & Terminal voltage of the battery & $\mathrm{V}$ & Algebraic \\
$i$ & Total current through the battery & $\mathrm{A}$ & Algebraic \\
$i_{\text{SEI}}$ & Current due to SEI side reaction & $\mathrm{A}$ & Algebraic \\
$\theta_p$, $\theta_n$ & Normalized surface lithium concentration ratio & & Algebraic \\
$U_p$, $U_n$ & Open-circuit potentials of positive/negative electrodes & $\mathrm{V}$ & Algebraic \\
$J_p$, $J_n$ & Reaction current density per unit electroactive surface area & $\mathrm{A}/\mathrm{m}^2$ & Algebraic \\
$D_p$, $D_n$ & Solid phase diffusivity (positive/negative electrodes) & $\mathrm{m}^2/\mathrm{s}$ & Constant \\
$R_p$, $R_n$ & Particles radius (positive/negative electrodes) & $\mathrm{m}$ & Constant \\
$\mathcal{F}$ & Faraday’s constant & $\mathrm{96,485 C/mol}$ & Constant \\
$a_p, a_n$ & Particle surface area to volume (positive/negative electrodes) & $\mathrm{m}^2/\mathrm{m}^3$ & Constant \\
$l_p, l_n$ & Region thickness (positive/negative electrodes) & $\mathrm{m}$ & Constant \\
$k_p, k_n$ & Rate for lithium intercalation reaction (positive/negative electrodes) & $\mathrm{m}^{2.5}/(\mathrm{mol}^{0.5} \cdot \mathrm{s})$ & Constant \\
$C_e$ & Electrolyte lithium-ion concentration & $\mathrm{mol}/\mathrm{m}^3$ & Constant \\
$\kappa_{\text{SEI}}$ & Ionic conductivity of SEI layer & $\mathrm{S/m}$& Constant \\
$T$ & Temporature & $\mathrm{K}$ & Constant \\
$\rho_{\text{SEI}}$ & Average density of the SEI layer & $\mathrm{Kg}/\mathrm{m}^3$ & Constant \\
$M_{\text{SEI}}$ & Molecular weight of the SEI layer & $\mathrm{Kg}/\mathrm{mol}$ & Constant \\
$U_{\text{ref}}$ & SEI growth onset voltage at negative electrode & $\mathrm{V}$ & Constant \\
$C_{p}^{\text{max}}$, $C_{n}^{\text{max}}$ & Maximum solid phase concentration (positive/negative electrodes) & $\mathrm{mol}/\mathrm{m}^3$ & Constant \\
\bottomrule
\end{tabularx}
\end{table*}

High-fidelity lithium-ion battery models are typically formulated using the SPM \cite{santhanagopalan2006review}. By neglecting electrolyte concentration gradients and simplifying radial diffusion in active particles, the SPM significantly reduces computational complexity while retaining the essential physics for tracking SOC and degradation. This simplification makes it particularly suitable for MPC, especially at moderate or low charge/discharge power levels where the simplifications introduce minimal error.

The SPM equations presented in this section are adapted from \cite{cao2020multiscale}, and the corresponding parameter set stems from a comprehensive identification campaign on ${LiFePO}_4$ cells via repeated voltage/current cycling experiments reported in~\cite{Forman2012}. ~\Cref{tab:sp_variables} summarizes the variables and parameters used in the SPM. The variables classified as state evolve according to differential equations, while those labeled as algebraic are determined solely by algebraic constraints. The input variables correspond to the control inputs of the system. All remaining quantities are fixed parameters,.

In SPM, assuming isothermal operation, the governing equations can be organized into five major components: 
lithium concentration dynamics within the electrode particles, the open-circuit 
voltages (OCVs) of the electrodes, interfacial electrochemical reaction kinetics, 
side reactions and Solid Electrolyte Interphase (SEI) film growth, and the 
resulting capacity fade and energy output. Each component collectively captures 
the interplay between diffusion, reaction thermodynamics, degradation, and 
macroscopic power, as detailed in the following key aspects.

\noindent\textbf{Lithium Concentration Dynamics}

\begin{align}
SOC &= \frac{C_{n}^{\text{avg}}}{C_{n}^{\text{max}}} \\
\frac{d C_{j}^{\text{avg}}}{dt} &= -15\,\frac{D_j}{R_{j}^2}\,(C_{j}^{\text{avg}} - C_{j}^{\text{surf}}), j\in\{p,n\}\label{eq:formala start} \\
C_{p}^{\text{surf}} - C_{p}^{\text{avg}}&=- \frac{R_{p}}{5}\times \frac{i}{\mathcal{F}\,D_p\,a_p\,l_p}\label{eq:concentrateDiffer positive}  \\
C_{n}^{\text{surf}} - C_{n}^{\text{avg}}&= \frac{R_{n}}{5}\times \frac{i_{int}}{\mathcal{F}\,D_n\,a_n\,l_n}\label{eq:concentrateDiffer negative} 
\end{align}
The temporal evolution of the average lithium concentration in the spherical particles of the positive ($p$) and negative ($n$) electrodes is described by ~\Cref{eq:formala start}. This equation represents the dominant diffusion-driven dynamics inside each particle. The relationships in ~\Cref{eq:concentrateDiffer positive,eq:concentrateDiffer negative} connect the average concentration to the surface concentration by imposing flux balance at the particle surface, thereby linking the diffusion field to the macroscopic current flow. These simplified expressions, which are actually derived as approximations from the boundary flux balance condition rather than by solving the full diffusion equation, eliminate the need to explicitly solve the full radial diffusion equation, while retaining the essential physics needed to track the electrodes’ SOC with adequate accuracy for control purposes.

\noindent\textbf{Open-Circuit Voltage (OCV)}
Based on empirical fits to OCV–SOC measurements for LiFePO\(_4\)-based cells, \(U_p\) and \(U_n\) are expressed as:
\begin{align}
\theta_j &= \frac{C^{\text{surf}}_{j}}{C_{j}^{\max}}, \quad j \in \{p,n\}, \\
U_j &= f_j(\theta_j), \quad j \in \{p,n\}, 
\end{align}
where \(f_p\) and \(f_n\) are smooth functions capturing the nonlinear dependence of the OCV on the state-of-charge. 

\noindent\textbf{Butler–Volmer Reaction Currents}

The intercalation and deintercalation of lithium ions at the electrode–electrolyte interfaces are modeled using Butler–Volmer kinetics. ~\Cref{eq:Butler-Volmer-positive,eq:Butler-Volmer-negative} express the reaction current densities at the positive and negative electrodes as nonlinear functions of the local surface concentrations, overpotentials, and reaction rate parameters. By linking these reaction currents with the concentration dynamics, the model establishes a consistent bridge between the applied current and the underlying electrochemical processes, which ultimately govern charge transfer and polarization behavior.
\begin{align}
\dfrac{i}{a_p\,\mathcal{F}\,l_p} &= 2\,k_p\,C_e^{0.5}\,(C_{p}^{\text{max}} - C_{p}^{\text{surf}})^{0.5}\,(C_{p}^{\text{surf}})^{0.5}\,
\sinh\!\biggl(0.5\,\dfrac{\mathcal{F}(\phi_p - U_p)}{R\,T}\biggr),
\label{eq:Butler-Volmer-positive} \\[0.5em]
-\dfrac{i_{\text{int}}}{a_n\,\mathcal{F}\,l_{n}} &= 2\,k_n\,C_e^{0.5}\,(C_{n}^{\text{max}} - C_{n}^{\text{surf}})^{0.5}\,(C_{n}^{\text{surf}})^{0.5}\,
\sinh\!\biggl(0.5\,\dfrac{\mathcal{F}}{R\,T}\Bigl(\phi_n - U_n + \dfrac{\delta_{\text{SEI}}}{\kappa_{\text{SEI}}}\dfrac{i}{a_n\,l_{n}}\Bigr)\biggr).
\label{eq:Butler-Volmer-negative}
\end{align}

\noindent\textbf{SEI Side Reaction, Film Growth and Capacity Fade}

At the negative electrode, parasitic side reactions lead to the formation and growth of the Solid Electrolyte Interphase (SEI) film. The SEI reaction current $i_{\text{SEI}}$ and the evolution of the film thickness $\delta_{\text{SEI}}$ are given by the governing equations above. These reactions not only consume active lithium, reducing the cell’s available capacity over time, but also increase interfacial resistance, altering the battery’s dynamic response. By integrating the SEI current over time, the model captures cumulative capacity fade $c_f$, allowing it to predict the gradual degradation of the battery over its operating life.
\begin{align}
i_{\text{SEI}} &= a_n\,l_{n}\,k_{\text{SEI}}\,
\exp\!\Biggl(-\dfrac{\mathcal{F}}{R\,T}\Bigl(\phi_n - U_{\text{ref}} + \dfrac{\delta_{\text{SEI}}}{\kappa_{\text{SEI}}}\,\dfrac{i}{a_n\,l_{n}}\Bigr)\Biggr),
\label{eq:sei_current} \\[0.5em]
\frac{d\delta_{\text{SEI}}}{dt} &= \frac{i_{\text{SEI}}\,M_{\text{SEI}}}{\mathcal{F}\,\rho_{\text{SEI}}\,a_n\,l_{n}},
\label{eq:sei_growth} \\[0.5em]
\frac{d c_f}{dt} &= \frac{i_{\text{SEI}}}{3600}.
\label{eq:cf_evolution}
\end{align}

\noindent\textbf{Current, Voltage and Power}

The terminal power output of the cell is determined by the product of the terminal voltage and the total applied current, as expressed in ~\Cref{eq:Power Conservation}. The total current $i$ is the sum of the intercalation current $i_{\text{int}}$ and the parasitic SEI side reaction current $i_{\text{SEI}}$, as given by the current balance relation:
\begin{align}
i &= i_{\text{int}} + i_{\text{SEI}}\\
P &= (\phi_{n} - \phi_{p})\times i \label{eq:Power Conservation}
\end{align}
The charging power $P$ is often described using the C-rate, which indicates how fast a battery is charged relative to its capacity—for example, 1C means full charge in one hour, 2C in half an hour.

\vspace{0.5em}

\noindent\textit{Collectively,} this formulation ensures current conservation while consistently linking the internal electrochemical processes—both reversible (intercalation) and irreversible (SEI growth)—to the macroscopic power flow at the cell terminals.

\subsection{Kriging Surrogate} \label{sec:Surrogate}
\Crefrange{eq:formala start}{eq:Power Conservation} form a Differential-Algebraic Equation (DAE) system and were implemented in earlier work by Cao et al.\cite{cao2020multiscale} for the arbitrage decision of a single battery. However, BSS involve multiple battery models and discrete decisions regarding which battery to swap, making even SPM still overly complex. In particular, the presence of exponential and hyperbolic sine terms often leads to severe numerical instability, making the problem extremely difficult to solve efficiently. To tackle this challenge, we adopt a surrogate modeling approach based on a Kriging model. This significantly simplifies the optimization problem while largely preserving the model's accuracy.

Kriging, originally developed in the field of spatial statistics, is a powerful spatial interpolation technique used for making optimal predictions at arbitrary spatial locations \cite{jones1998efficient,sacks1989design}.
Due to its high prediction accuracy and low number of function evaluations \cite{lin2017process}, Kriging is well-suited to serve as a surrogate model. Its underlying model is:
\begin{align}
f(x) &= \mu(x) + \delta(x)
\end{align}
Here, $\mu(x)$ represents the deterministic trend component (mean function), and $\delta(x)$ is a zero-mean stochastic process that models spatial correlation through a covariance function. In practice, the mean function is often chosen as a constant or a linear function. In this work, we adopt a constant mean function, i.e., $\mu(x) = \mu_0$. The covariance function, also known as the kernel, encodes the similarity between input points and plays a key role in defining the smoothness and generalization ability of the surrogate model. 

Given a set of training inputs $X$ and corresponding outputs $y=f(X)$, the Kriging predictor at a new input $x_*$ follows a normal distribution with the following mean and variance:
\begin{align}
\mu_* &= \mu + k_* K^{-1} \left(y - \mu 1\right) \\
\sigma_* &= k_{**} - k_* K^{-1} k_*^\top
\end{align}
In these expressions:
\begin{itemize}
\item $K$ is the covariance matrix of the training samples,
\item $k_*$ is the covariance vector between the test point $x_*$ and the training inputs $X$,
\item $k_{**}$ is the self-covariance of the test point.
\end{itemize}
All covariance terms are computed using the selected kernel function $k(\cdot, \cdot)$. A commonly used kernel is the Radial Basis Function (RBF), also known as the squared exponential kernel, defined as:
\begin{align}
k(x,x') = \prod_{d=1}^{D} \exp\,\!\bigl(-\theta_{d}|x_{d}-x'_{d}|^{2}\bigr),
\end{align}
where:
\begin{itemize}
\item $x, x'$ are input vectors, and D is the dimension of $x$.
\item $\theta_{d}$ is the length scale parameter, which controls the sensitivity of the function to differences in the input;
\end{itemize}

\noindent\textbf{Kriging Surrogate for SPM}

For the battery model, the state vector $x$ consists of the four key state variables $C^{\text{avg}}_p$, $C^{\text{avg}}_n$, $\delta_{\text{SEI}}$, and $c_f$, while the control input $u$ corresponds to the power input $P$. By retaining only these dominant dynamic states and omitting numerous algebraic variables, we simplify the original DAE system into a tractable ordinary differential equation (ODE) model, which serves as the basis for the surrogate construction. The Kriging surrogate is trained to approximate the state increment, formulated as
\begin{align}
\Delta x_{i+1} &= f(x_{i+1}, u_{i+1}) \label{eq:state transition 1}
\end{align}
where $\Delta x_{i+1}$ represents the change in the state between two consecutive time steps i and i+1.
This incremental formulation is particularly beneficial for accurately capturing the monotonically increasing behavior of $\delta_{\text{SEI}}$ and $c_f$, as directly modeling the absolute states tends to introduce bias and reduce fidelity for such slowly varying, cumulative quantities.
The actual state evolution is then obtained by integrating this increment:
\begin{align}
x_{i+1} &= x_i + \Delta x_{i+1} \label{eq:state transition 2}
\end{align}
It is worth noting that we fit the state transition model in the form of backward integration. This choice leads to a more compact representation of the entire BSS model, as will be further elaborated in ~\Cref{sec:Multiscale_Optimization_of_BSS}.

\begin{figure}[!h]
    \centering
    \begin{subfigure}[b]{0.49\linewidth}
        \includegraphics[width=\linewidth]{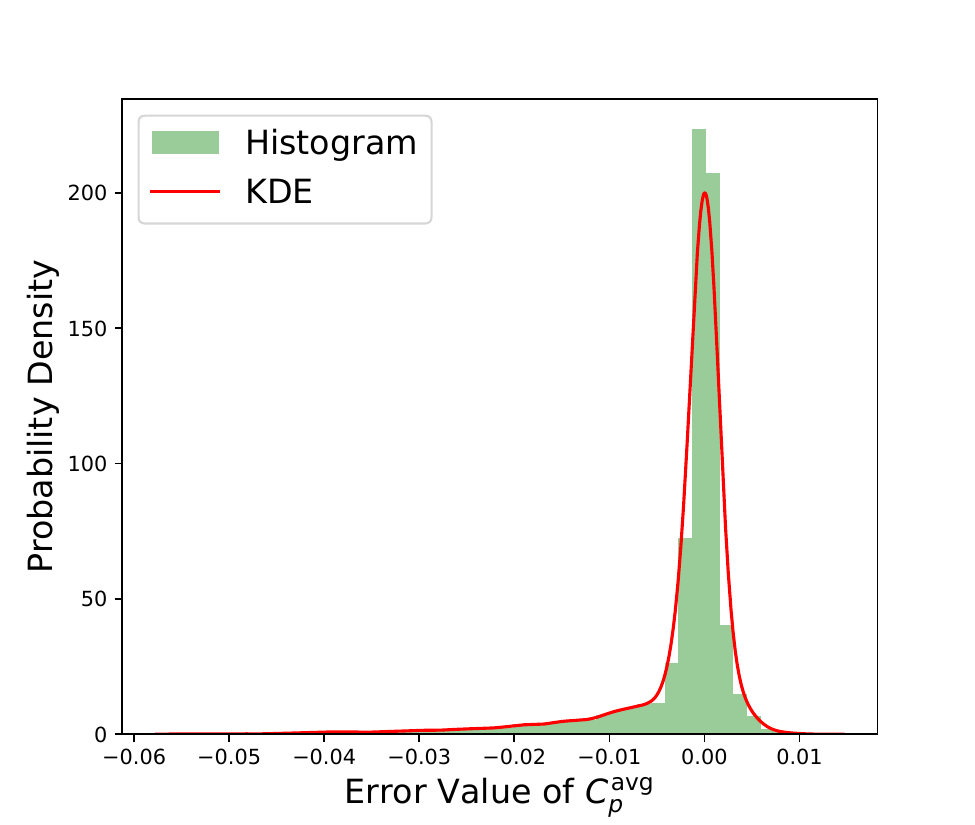}
        \caption{$C^{\text{avg}}_p$ fitting residual.}
        \label{fig:relative_error_csp}
    \end{subfigure}
    \begin{subfigure}[b]{0.49\linewidth}
        \includegraphics[width=\linewidth]{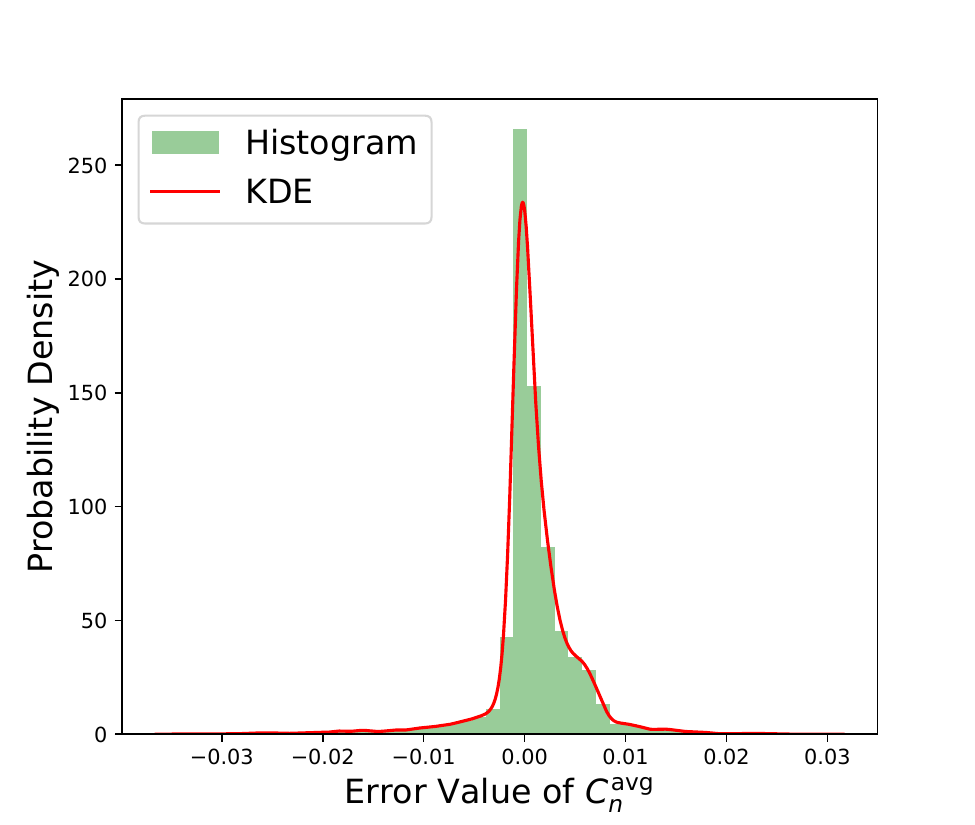}
        \caption{$C^{\text{avg}}_n$ fitting residual.}
        \label{fig:relative_error_csn}
    \end{subfigure}
    \begin{subfigure}[b]{0.49\linewidth}
        \includegraphics[width=\linewidth]{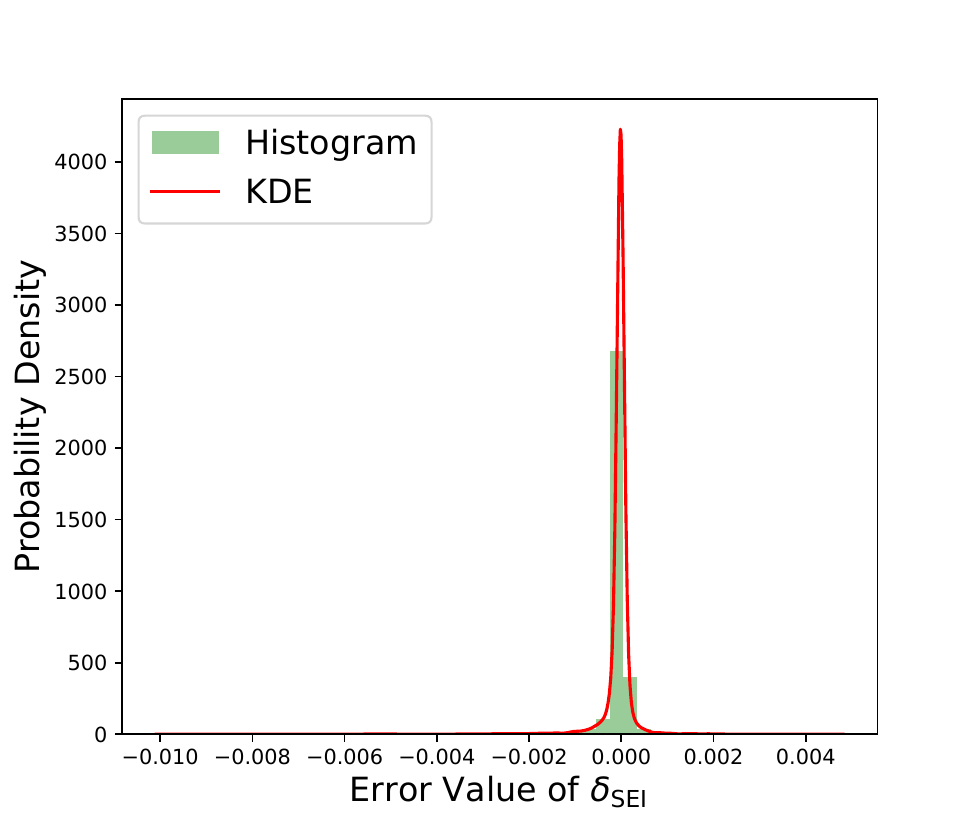}
        \caption{$\delta_{\text{SEI}}$ fitting residual.}
        \label{fig:relative_error_delta_sei}
    \end{subfigure}
    \begin{subfigure}[b]{0.49\linewidth}
        \includegraphics[width=\linewidth]{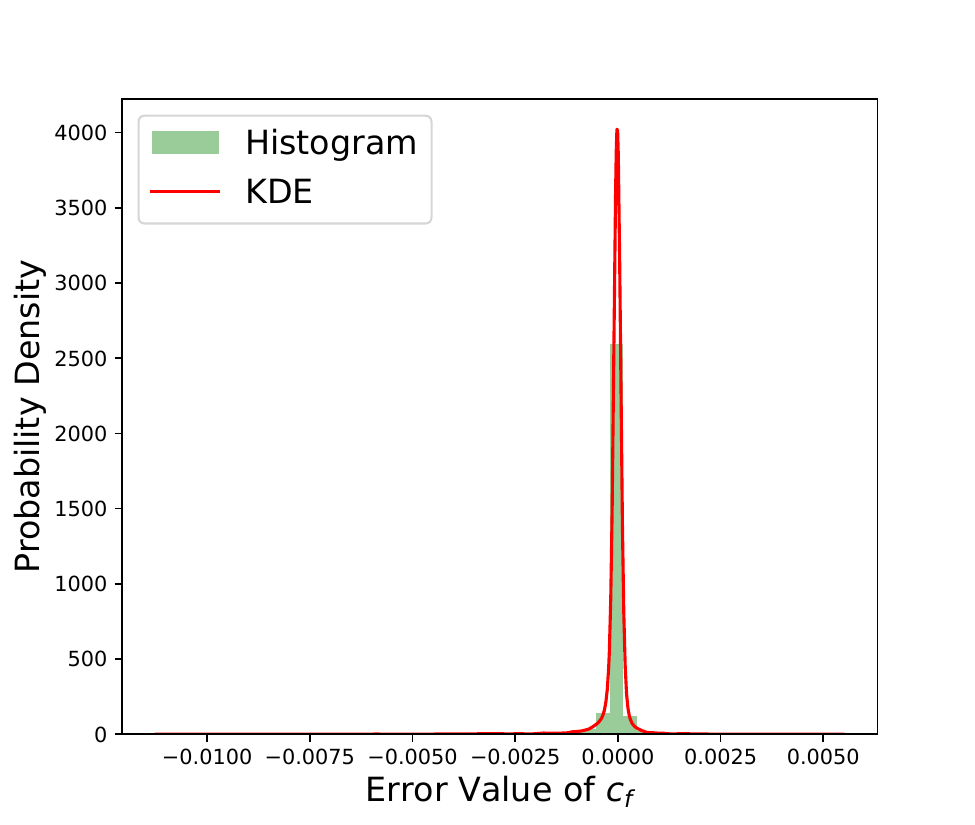}
        \caption{$c_f$ fitting residual.}
        \label{fig:relative_error_cf}
    \end{subfigure}
    \caption{Histogram of Kriging fitting for $C^{\text{avg}}_p$, $C^{\text{avg}}_n$, $\delta_{\text{SEI}}$, and $c_f$}
    \label{fig:relative_error}
\end{figure}

Input-output data pairs are generated offline by directly performing numerical integration on the DAE system ~\Crefrange{eq:formala start}{eq:Power Conservation} using the IDA solver\cite{rackauckas2017differentialequations,gardner2022sundials,hindmarsh2005sundials}. We constructed a dataset by numerically integrating the evolution of a battery from its initial state through randomized charge–discharge cycles until reaching end-of-life. Based on this high-fidelity dataset, a surrogate model is trained to approximate the battery dynamics. By omitting the algebraic variables present in the original electrochemical model, we reduce the system to a set of ODEs, lowering the computational complexity while maintaining accuracy. Importantly, since the surrogate model is trained on data sampled from the original continuous-time model, its prediction error remains consistent with the true dynamics regardless of the discretization step size. This allows us to use a relatively large time step when solving optimal control problems without sacrificing accuracy.

We use maximum likelihood estimation (MLE) to train the kriging model. The fitting accuracy of the resulting Kriging surrogate model is illustrated in ~\Cref{fig:relative_error}, where the relative error distribution is visualized using kernel density estimation. As shown, the surrogate achieves excellent agreement with the original model across all state variables. In particular, the slowly varying degradation-related states $\delta_{\text{SEI}}$ (\Cref{fig:relative_error_delta_sei}) and $c_f$ (\Cref{fig:relative_error_cf}) exhibit relative errors concentrated within $\pm 0.002$, indicating highly accurate approximation. For the more dynamic concentration-related states $C^{\text{avg}}_p$ (\Cref{fig:relative_error_csp}) and $C^{\text{avg}}_n$ (\Cref{fig:relative_error_csn}), the surrogate maintains errors within $\pm 0.03$, with the error distribution sharply centered around zero and no noticeable bias. 

While the MPC framework requires repeated trajectory rollouts over the prediction horizon, only the first control input is applied at each time step during closed-loop execution. This receding-horizon structure inherently mitigates the impact of small prediction errors. Overall, the Kriging surrogate offers a computationally efficient yet sufficiently accurate representation of the original SPM, capturing both state-of-charge dynamics and capacity fade, and thus providing a reliable foundation for multiscale BSS operational optimization.

\subsection{Optimal Control of BSS} \label{sec:Multiscale_Optimization_of_BSS}

To account for both fast operational dynamics and slow battery aging, we formulate a multiscale model of the BSS that integrates swapping logistics, electrochemical behavior, and degradation. Viewing this as a control system, we pose a finite-horizon mixed-integer optimal control problem that links short-term decisions with long-term health, forming the basis of the predictive strategy in the next section.

\subsubsection{State Transition Model} \label{sec:State Transition Model}

We consider a deterministic framework for modeling the dynamics of a BSS, where exogenous factors — such as real-time electricity prices and vehicle arrival patterns — are assumed to be known a priori and treated as time-dependent inputs to the system. The analysis is based on two primary datasets: historical Day-Ahead Market energy prices from PJM Interconnection for a one-year period and historical hourly customer swapping demand records from NIO's BSS located at the Shanghai Changtai International Financial Center. Furthermore, we assume that the batteries used in EVs discharge in a stable manner within their appropriate operating range, and thus capacity degradation is neglected. To capture the hourly fluctuations in electricity prices and to simplify computation, we discretize the system dynamics using a 1-hour time step ($\Delta t = 1 \text{h}$). A 1-hour interval is sufficient for full battery charging under the considered C-rate conditions. The state and control variables at these grid points are discretized as follow,
\begin{align*}
    x_i &= x(t_i), \quad i = 0, \dots, N, \\
    u_i &= u(t_i), \quad i = 1, \dots, N.
\end{align*}

Let $K$ be the total number of batteries in the station. The state vector at time step $i$ is defined as:
\begin{align*}
x_i =
\begin{bmatrix}
x_{1,i}^\top & x_{2,i}^\top & \cdots & x_{K,i}^\top
\end{bmatrix}^\top \in \mathbb{R}^{4\times K},
\end{align*}
where the state of each battery $k$ is:
\begin{align*}
x_{k,i} =
\begin{bmatrix}
C^{\text{avg}}_{p,k,i}, \quad
C^{\text{avg}}_{n,k,i}, \quad
\delta_{\text{SEI},k,i}, \quad
c_{f,k,i}
\end{bmatrix}^\top .
\end{align*}

The control input at time step $i$ is:
\begin{align*}
u_i =
\begin{bmatrix}
P_{1,i} & P_{2,i} & \cdots & P_{K,i}
\end{bmatrix}^\top \in \mathbb{R}^K,
\end{align*}
where $P_{k,i}$ denotes the net power exchanged between battery $k$ and the grid at step $i$ (positive for discharging/selling, negative for charging/buying).

To model swapping operations, we introduce a binary selection vector for step 1 to N:
\begin{align*}
b_i =
\begin{bmatrix}
b_{1,i} & b_{2,i} & \cdots & b_{K,i}
\end{bmatrix}^\top \in \{0,1\}^K,
\end{align*}
where $b_{k,i} = 1$ if battery $k$ is selected for swapping at step $i$, and $b_{k,i} = 0$ otherwise. The number of batteries to be swapped, denoted $S_i$, is assumed to be known at each step $i$, with the constraint:
\begin{align}
\sum_{k=1}^K b_{k,i} = S_i.
\end{align}

The discrete-time state transition with swapping is expressed as:
\begin{align}
x_{i+1} = b_{i} \cdot x_{\text{swap}} + (1 - b_{i}) \cdot x_i + \Delta x_{i+1}, \label{eq:jump_state_model}
\end{align}
where:
\begin{itemize}
    \item $x_{\text{swap}}$ is the state of a freshly swapped battery, assumed to be fixed and known. 
    \item $\Delta x_{i+1}$ represents the state increment from system dynamics, as described in ~\Crefrange{eq:state transition 1}{eq:state transition 2}.
\end{itemize}

Using a surrogate model to approximate the backward state update (via $\Delta x_{i+1}$) instead of a forward one simplifies the formulation. If the forward form were used, we would need to track two separate increments, $\Delta x_{\text{swap}}$ and $\Delta x_i$, leading to:
\begin{align*}
\Delta x_i &= f(x_i, u_i), \\
x_{i+1} &= (1 - b_i) \cdot (x_i + \Delta x_i) + b_i \cdot (x_{\text{swap}} + \Delta x_{\text{swap}}),
\end{align*}
which significantly complicates the model.

To maintain physical realism, only batteries with an SOC above a minimum threshold are eligible for swapping:
\begin{align}
SOC_i \geq b_i \cdot (\tilde{SOC} + \varepsilon),
\end{align}
ensuring that swapping occurs only when a battery has sufficient charge. The term \(\varepsilon\) compensates for the Kriging-induced approximation error and ensures the feasibility of swapping constraints during numerical optimization.

\subsubsection{Objective Function} \label{sec:Objective Function}

The optimization objective balances economic profitability with long-term battery health preservation. 
It combines a reward term $\mathcal{R}$ and a penalty term $\mathcal{P}$, forming the total objective:
\begin{align}
\Phi = \mathcal{R} - \mathcal{P},
\end{align}
which we seek to maximize over the discrete planning horizon $i = 0, \dots, N$.

\paragraph{Reward Term.}
The reward function encourages economically favorable grid interactions, where selling electricity to the grid is rewarded and energy purchases are penalized. It is defined as:
\begin{align}
\mathcal{R} = \sum_{k=1}^K \sum_{i=1}^{N}  P_{k,i} \cdot \rho_i,
\end{align}
where $\rho_i$ is the real-time price at step $i$.

\paragraph{Penalty Term.}
The penalty function discourages behaviors that lead to accelerated battery degradation and imbalanced utilization across the battery fleet. It consists of two components:
\begin{align}
\mathcal{P} &= w_1 \mathcal{P}_1 + w_2 \mathcal{P}_2,
\end{align}
where $w_1$ and $w_2$ are user-defined weights.

\begin{itemize}
    \item \textbf{Capacity Fade Penalty:}  
    This term penalizes cumulative capacity loss, modeled through the increase in the degradation variable $c_{f,k}(t)$:
    \begin{align}
        \mathcal{P}_1 &= \sum_{k=1}^K \Pi \cdot \bigg((1 - b_{k,i}) \cdot\big(c_{f,k,i+1} - c_{f,k,i}\big) + b_{k,i} \cdot\big(c_{f,k,i+1} - c_{f,swap}\big)\bigg),
    \end{align}
    where $\Pi$ is the economic cost per unit capacity fade, and $c_{f,swap}$ is the capacity fade of freshly swapped battery.

    \item \textbf{Utilization Balance Penalty:}  
     Due to the nonlinear aging behavior of lithium-ion batteries, repeatedly using a small subset can lead to disproportionate degradation. To address this, batteries with higher capacity fade are preferentially selected for swapping, encouraging uniform degradation across the entire fleet. We define:
    \begin{align}
        \mathcal{P}_2 &= \sum_{k=1}^K (1 - b_{k,i}) \cdot \big( c_{f,k,i} - c_{f,\mathrm{offset}} \big),\label{eq:distribution penalty}
    \end{align}
     where $c_{f,\text{offset}}$ represents a fixed constant used to normalize degradation levels, typically chosen as $\min_k c_{f,k}(0)$.
\end{itemize}

By jointly maximizing economic returns and extending system lifespan, this formulation supports both profitability and sustainability in BSS operations.

\subsubsection{Discrete-Time Mixed-Integer Optimal Control Problem} \label{sec:Discrete_OCP}

The above formulation captures a multiscale optimization problem defined over a finite time horizon. To solve this, we formulate the BSS operation problem as a finite-horizon mixed-integer optimal control problem (MIOCP).

For notational compactness, we first note that the discrete-time system dynamics ~\Cref{eq:jump_state_model} can be interpreted as continuity constraints linking consecutive nodes in the horizon. Specifically, these continuity constraints are written as
$x_i = F_i(x_{i+1}, u_{i+1}, b_{i})$,
and the stage cost at each node is denoted as $L_i = L_i(x_i,u_i,b_i)$, where $i = 0,\dots,N-1$.

This leads to the following finite-horizon MIOCP formulation:

\begin{align}
    \max_{\{x_i,u_i,b_i\}} \quad & \sum_{i=1}^{N} \Phi_i(x_i, u_i, b_i) \label{eq:surrogate_ocp_start} \\
\text{s.t.} \quad & x_i = F_i(x_{i+1}, u_{i+1}, b_{i}), \\
& x_0 = x_{\text{init}}, \\
& SOC_i \geq b_{i} \cdot (\tilde{SOC} + \varepsilon), \\
& b_i \in \{0,1\}^{K}, \\
& \sum_{k=1}^{K} b_{k,i} = S_i, \quad i = 1,\dots,N,\label{eq:surrogate_ocp_end} \\
& \text{other path, state, and control constraints.}
\end{align}

The decision variables of the optimal control are all discretized states, inputs, and binary decisions:
\begin{align*}
\{x_0, x_1, \dots, x_N, \; u_1, \dots, u_N, \; b_1, \dots, b_N\}.
\end{align*}

\subsection{BSS-MPC} \label{sec:BSS-MPC}

To enable practical and real-time solution of the finite-horizon mixed-integer optimal control problem described above, we adopt a MPC strategy. In our implementation, a prediction horizon of 24 hours (one full day) is selected to fully account for the fluctuations in the day-ahead electricity market as well as the variations in battery swapping demand. Rather than optimizing over the entire scheduling period at once, MPC solves the problem over this moving horizon and updates the solution at each time step, thereby adapting online to changes in system states and external signals.

According to multiple shooting method~\cite{bock1984multiple}, in each MPC iteration, the finite-horizon optimal control problem is directly transcribed into a MINLP problem, effectively capturing the coupling between discrete decisions and continuous system dynamics within a unified framework.

In practice, the bilinear coupling between binary and continuous variables in ~\Cref{eq:jump_state_model} poses challenges for numerical optimization. To enable tractable formulation, we apply the Big-M method to obtain an equivalent linear representation:
\begin{align}
x_{i+1} - \Delta x_{i+1} &\leq x_i + M \cdot b_{i}, \\
x_{i+1} - \Delta x_{i+1} &\geq x_i - M \cdot b_{i}, \\
x_{i+1} - \Delta x_{i+1} &\leq x_{\text{swap}} + M \cdot (1 - b_{i}), \\
x_{i+1} - \Delta x_{i+1} &\geq x_{\text{swap}} - M \cdot (1 - b_{i}),
\end{align}
where $M$ is a sufficiently large constant, typically chosen as $M = |x_{\max} - x_{\min}|$. These constraints are logically equivalent to the following switching behavior:
\begin{align*}
x_{i+1} =
\begin{cases}
x_i + \Delta x_{i+1}, & b_{i} = 0 \quad \text{(no swap)}, \\
x_{\text{swap}} + \Delta x_{i+1}, & b_{i} = 1 \quad \text{(swap occurs)}.
\end{cases}
\end{align*}

\RestyleAlgo{ruled}
\begin{algorithm}[!hbt]
\caption{BSS-MPC} \label{algo:BSS-MPC}
\SetAlgoLined
\KwData{Initial grid price $\rho$, swap demand $S$, initial state $x_0$}
\KwResult{Control input sequence $u$ and state trajectory $x$}

\BlankLine
\textbf{Step 1: Kriging Surrogate Fitting}\;
Sample training state points $X$ from training set $\mathscr{X}_{\text{train}}$\;
\While{relative error (RE) $\geq$ Tol}{
    Add new points $X_*$ and update Kriging parameter $\theta_d$\;
    Compute RE on testing set $\mathscr{X}_{\text{test}}$, $RE \gets \frac{\lVert y - \tilde{y} \rVert}{\lVert y \rVert}$\;
}

\BlankLine
\textbf{Step 2: MPC Loop}\;
$i \gets 0$\;
\While{$i <$ total simulation time}{
    Define prediction horizon as $[i, i+N]$\;
    Solve surrogate OCP ~\Crefrange{eq:surrogate_ocp_start}{eq:surrogate_ocp_end} over $[i, i+N]$\;
    Perform feasibility check using SPM ~\Crefrange{eq:formala start}{eq:Power Conservation}\;
    \If{solution not feasible}{
        Add new points $X_*$ and update Kriging parameter $\theta_d$\;
        \textbf{continue}\;
    }
    Apply the first control input $u_i$ to the real BSS system\;
    Obtain the resulting state $x_{i+1}$\;
    Update grid price $\rho_{i:i+N}$ and swap demand $S_{i:i+N}$\;
    $i \gets i + 1$\;
}
\end{algorithm}

By decoupling the state trajectories across intervals and treating integer variables explicitly, the formulation avoids severe nonlinearities and large-scale complexity, enabling the use of state-of-the-art MINLP solvers to handle the hybrid nature of the problem more efficiently. The MINLP model is implemented using the JuMP package in the Julia language~\cite{Lubin2023} and solved with the DICOPT solver~\cite{grossmann2002gams}. We configured DICOPT to use IPOPT~\cite{wachter2006implementation} as the NLP solver and CPLEX~\cite{manual1987ibm} as the MILP solver.

The BSS-MPC framework, as outlined in ~\Cref{algo:BSS-MPC}, comprises two main stages: offline surrogate model construction and online receding-horizon control. In the offline stage, a Kriging surrogate is trained using an initial set of sampled data points and iteratively refined via cross-validation until the prediction error falls below a prescribed tolerance. During the online stage, at each control step, an OCP is formulated over the prediction horizon $[i, i+N]$ based on the prevailing electricity price and swapping demand. The solution is subsequently validated against the original model to ensure feasibility. If infeasibility is detected, additional data points are adaptively sampled in the vicinity of the current operating point, the surrogate is updated, and the OCP is re-solved. Once a feasible solution is obtained, only the first control action $u_i$ is implemented following the receding horizon strategy, the system state is updated through simulation (or in the real environment), and the procedure advances to the next time step.

\section{Experiment}\label{sec:Experiment}
To simulate a realistic operational scenario, we construct a virtual community with 200 batteries served by a single station. This configuration, with 179 EV drivers and 21 batteries held in reserve at the BSS, mimics a real-world scenario where a single BSS provides service for a neighborhood. The station maintains a constant inventory of 21 batteries, while the remaining batteries are in circulation with vehicle owners. Each battery is uniquely indexed from 1 to 200. Initially, batteries 1 through 21 are stored in the station, and batteries 22 through 200 are assigned to vehicles. The battery model used in the experiments is identical to that in \cite{cao2020multiscale}: A123 Systems ANR26650M1 cylindrical cells with $LiFePO_4$ cathodes are adopted, and all cell-level parameters have been thoroughly identified via extensive voltage/current cycling tests as reported in \cite{Forman2012}.

Battery swapping demands are handled in a first-in-first-out (FIFO) manner. When a vehicle arrives to swap a battery, it returns its current battery to the station and receives the next available one from the inventory. For example, if three swaps occur at the first time step, batteries 22 to 24 are returned to the station. The station, in turn, releases three batteries (e.g., batteries 4, 7, and 10), which are appended to the end of the vehicle battery queue. As a result, the new sequence of vehicle-held batteries becomes 25 to 200, followed by 4, 7, and 10.

Since battery swapping is a mandatory service, failing to provide a sufficiently charged battery could, in the long term, damage brand reputation and result in customer attrition. Therefore, while successful swapping does not generate direct revenue, failed swaps incur a cost. Specifically, we impose a penalty of \$1 for every 10\% shortfall in SOC relative to the swapping threshold. As in ~\Cref{sec:Multiscale_Optimization_of_BSS}, when a battery is used by an EV driver, it operates under stable discharge conditions, and thus capacity degradation is assumed negligible.

Moreover, because there is a mismatch between the simplified model used for planning and the high-fidelity DAE simulation used for evaluation, charge and discharge actions may sometimes result in overcharging or overdischarging. To mitigate such risks, protective strategies should be implemented to automatically halt charging or discharging once predefined limits are exceeded, thereby ensuring the system consistently operates within safe conditions.

\begin{figure}[!h]
    \centering

    \begin{subfigure}[b]{0.49\linewidth}
        \centering
        \includegraphics[width=\linewidth]{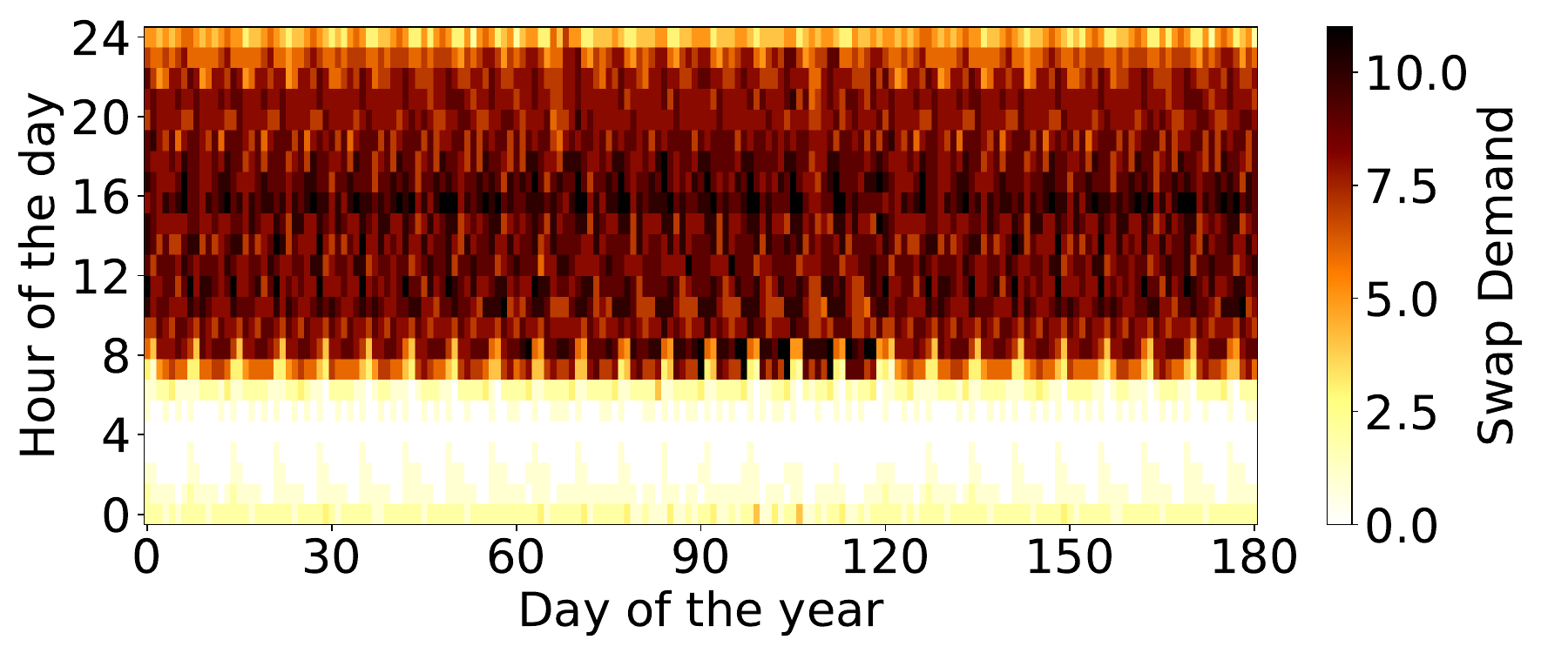}
        \caption{Swap demand by hour and day}
        \label{fig:heat_demand}
    \end{subfigure}
    \hfill
    \begin{subfigure}[b]{0.49\linewidth}
        \centering
        \includegraphics[width=\linewidth]{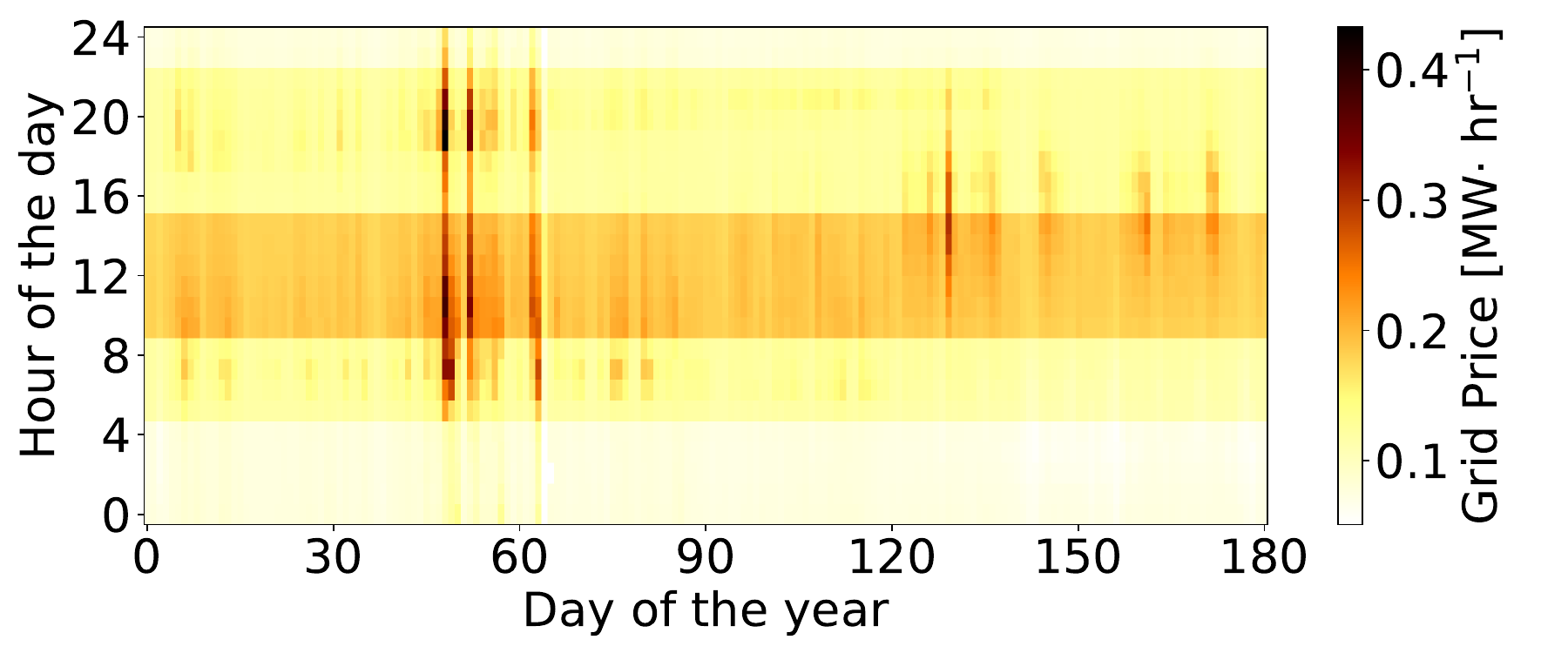}
        \caption{DAM electricity price by hour and day}
        \label{fig:heat_price}
    \end{subfigure}

    \caption{Temporal heatmaps of swap demand and electricity price: swap demand is concentrated during daytime and minimal overnight. Electricity prices exhibit a diurnal pattern with late-afternoon peaks, occasional spikes, and a mid-year volatility patch.}
    \label{fig:heatmap}
\end{figure}

The known variations in electricity prices and battery swapping demand are shown in ~\Cref{fig:heatmap}. It can be observed that swapping demand is generally higher during the daytime, while it decreases from evening to early morning; meanwhile, the peak electricity prices also occur during daytime working hours.


We additionally implement two models to highlight the advantages of our proposed approach in terms of economic efficiency and battery longevity:

\noindent\textbf{Low-Fidelity model}: This model simplifies battery-related states to a single variable—the SOC. Charging and discharging efficiencies are treated as constants, and capacity degradation is modeled as being proportional to the square of the power used. Rewards are assigned based on economic profit, while penalties are applied for capacity loss. This formulation reduces the optimal control problem to a mixed-integer quadratic programming (MIQP) problem:
\begin{equation}
\begin{split}
    \max_{\{x_i, u_i\}} \quad & \sum_{i=1}^{N} P_{i} \cdot \rho_i - \sum_{i=1}^{N} w \cdot P_i^2 \\
    \text{s.t.} \quad & SOC_{i+1} = b_{i} \cdot SOC_i + (1 - b_{i}) \cdot SOC_{\text{swap}} + \eta_{\text{ch}} \cdot P_{i+1},\\
    & SOC_i \geq b_{i} \cdot (\tilde{SOC} + \varepsilon) \\
    & b_i \in \{0,1\}^{K} \\
    & \sum_{k=1}^{K} b_{k,i} = S, \quad i = 1,\dots,N
\end{split}
\end{equation}
where \(w\) is the weight and \(\eta_{\text{ch}}\) is the charge and discharge efficiency. To ensure full charging within one discrete time step (1 hour), we set \(\eta_{\text{ch}} = \frac{1}{1C}\).

\noindent\textbf{Rule-Based model}: This model follows a heuristic strategy for battery dispatch and charging. Considering the nonlinear degradation characteristics of lithium-ion batteries, the strategy always selects the batteries with greater capacity loss for swapping. This helps to balance the overall aging across the community battery fleet, avoiding the accelerated degradation caused by overusing a subset of batteries and thereby extending the average lifespan of the entire system. Additionally, every battery swapped into the station is fully charged within one hour, thereby maintaining all batteries inside the station at full charge. Instead of using maximum-power charging, 
at each time step $i \in \{1, \ldots, N-1\}$, we solve the following root-finding problem to determine the control policy $u_{i+1}$:
\begin{equation}
\begin{split}
    \min_{\{u_i\}} \quad & \Vert x_{i} - F(x_{i+1},u_{i+1}) \Vert_2^2 \\
    \text{s.t.}\quad & SOC_{i+1} =\tilde{SOC} + \varepsilon
\end{split}
\end{equation}
where where $F(x_{i+1}, u_{i+1})$ is the backward Euler integration of SPM. Notably, we reformulate the original root-finding task as an optimization problem: the objective function is designed to drive the residual of the equation $ x_{i} - F(x_{i+1}, u_{i+1}) = 0 $ towards zero. In addition, unlike the BSS-MPC and low-fidelity strategies, the rule-based strategy merely considers a single time step rather than a full prediction horizon, which makes it the least computationally demanding approach.




The rule-based strategy ignores both capacity degradation and electricity price fluctuations. The low-fidelity method accounts for price variations but relies on a simplified charge–discharge model and penalizes power quadratically to approximate degradation. In contrast, BSS-MPC explicitly incorporates the aging model and surrogate state transitions, while simultaneously optimizing against dynamic electricity prices.

\subsection{Planning Results: A One-Day Case Study} \label{sec:planning_result}
\begin{figure*}
    \centering
    \includegraphics[width=0.9\linewidth]{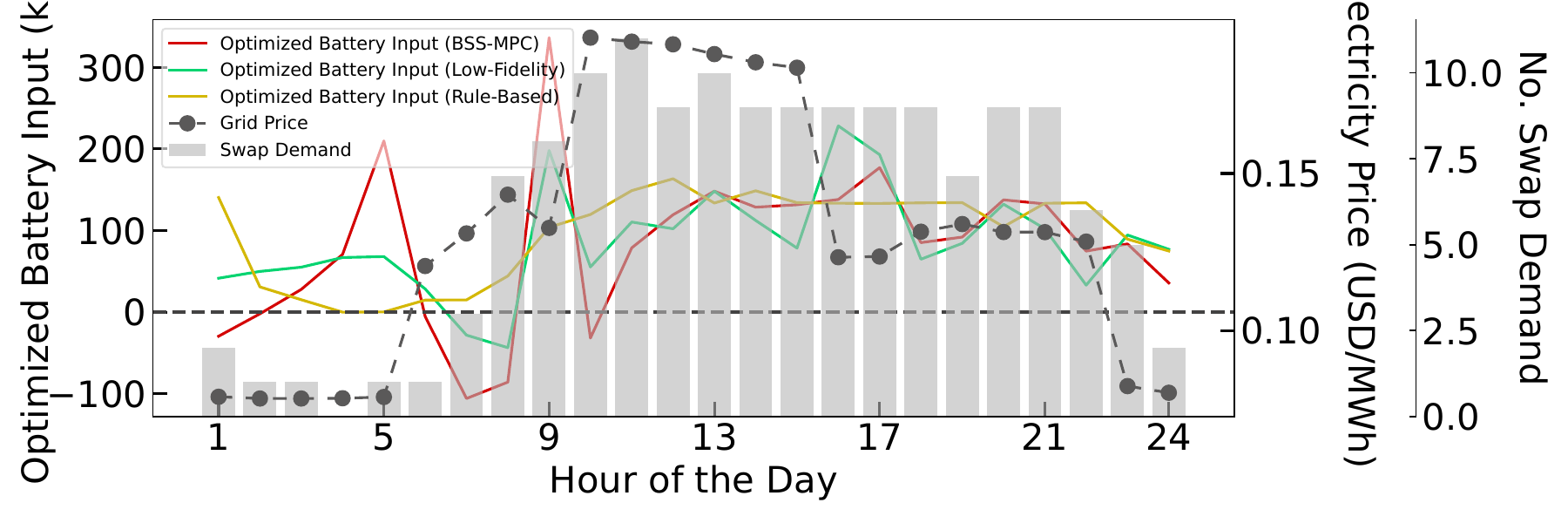}
    \caption{Within-horizon scheduling: optimized battery operation (MPC v.s. baselines) with grid price and swap demand}
    \label{fig:demand_price_optimized}
\end{figure*}
To ensure generality, we randomly selected a day (24 hours) to evaluate the effectiveness of the BSS-MPC planning.
~\Cref{fig:demand_price_optimized} presents the one-day optimization results of the BSS-MPC (Ours), Low-Fidelity, and Rule-Based methods, illustrating how the three approaches respond to fluctuations in electricity prices and swapping demand. As is evident, the Rule-Based method completely ignores price variations and directly charges the battery to $\tilde{SOC}$ after each swap, resulting in no arbitrage behavior.
In contrast, both the BSS-MPC and Low-Fidelity methods not only account for changes in demand but also adjust their charging and discharging strategies according to price fluctuations. During the first six hours, when both electricity prices and demand are low, the two methods choose to accumulate energy. This allows them to sell a certain amount of electricity for arbitrage during hours 7 and 8, when both prices and demand suddenly increase. At hour 9, there is a brief price dip, during which both methods charge significantly to prepare for the upcoming peak in price and demand. In particular, at hour 10, when the electricity price reaches its highest point, BSS-MPC is still able to sell some electricity.
From hour 11 onward, as swapping demand remains consistently high, the arbitrage opportunities shrink, and the decisions made by the three methods become similar. Since the price and demand peaks occur simultaneously, the optimal strategy is clear: the optimizer takes advantage of low-price, low-demand periods to gain small profits while purchasing enough electricity in advance to keep the battery at a high state of charge. This way, during peak periods, only a small amount of electricity needs to be purchased to meet the swapping demand, significantly reducing overall costs.

\begin{figure}[!h]
    \centering
    \begin{subfigure}[b]{0.7\linewidth}
        \centering
        \includegraphics[width=\linewidth]{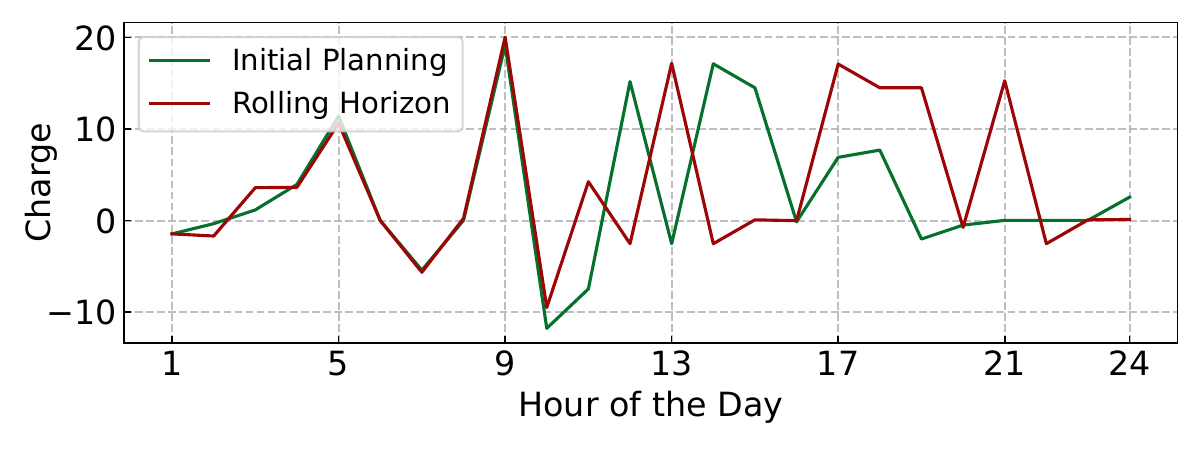}
        \caption{Charging Input}
        \label{fig:compare_input}
    \end{subfigure}

    \begin{subfigure}[b]{0.7\linewidth}
        \centering
        \includegraphics[width=\linewidth]{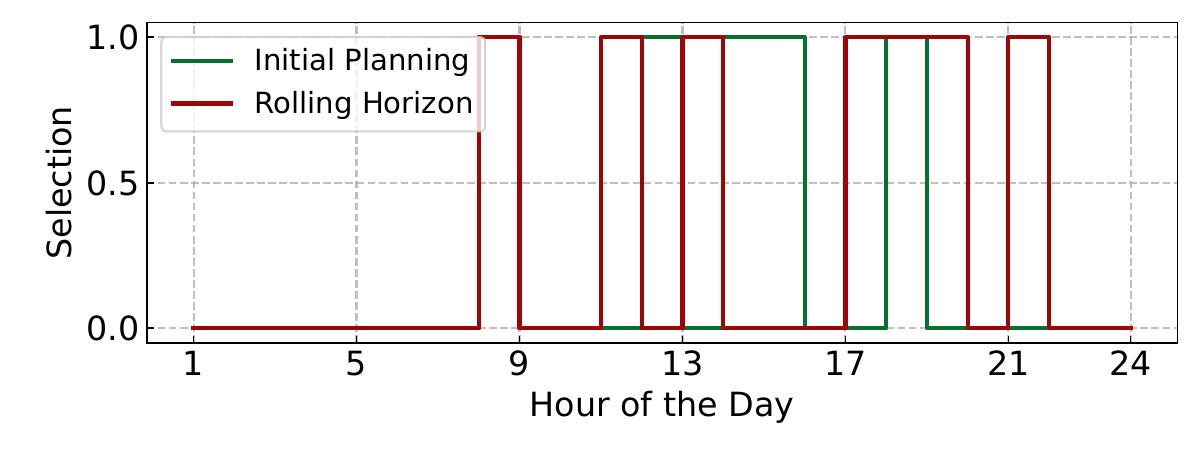}
        \caption{Selection}
        \label{fig:compare_selection}
    \end{subfigure}

    \begin{subfigure}[b]{0.7\linewidth}
        \centering
        \includegraphics[width=\linewidth]{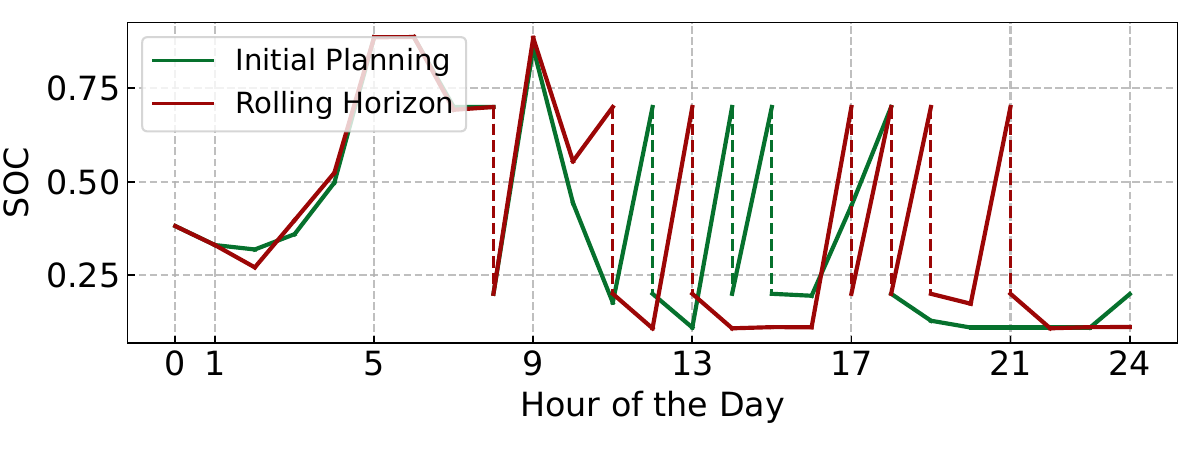}
        \caption{SOC}
        \label{fig:compare_soc}
    \end{subfigure}

    \caption{Plan vs. Realization (rolling horizon): charging power, swap decisions, and SOC over a day.}
    \label{fig:result_comparison}
\end{figure}

Additionally, we compared the results of the initial planning phase with those from the rolling optimization that covers the entire 24-hour period. Given that the battery swap station maintains 21 battery slots, displaying the results for all slots would be unclear. Therefore, we randomly selected one slot for analysis, as shown in~\Cref{fig:result_comparison}. ~\Cref{fig:compare_input} illustrates the charge and discharge input of the selected battery slot, while ~\Cref{fig:compare_selection} shows whether the battery from the selected slot is swapped out (with 1 indicating a swap). ~\Cref{fig:compare_soc} displays the evolution of the battery SOC in the selected slot, where the dashed lines indicate that the battery is swapped out, marking a state jump within the battery slot due to the swap.

It is evident that, in the initial hours (e.g., up to hour 10), the decisions made by the rolling optimization are quite similar to those from the initial planning phase. The battery selection decisions for this particular slot remain identical, and both the charge and discharge decisions follow similar trends with only minor adjustments. As a result, the SOC evolution for this slot also follows a similar path, as the underlying control decisions are aligned. This is because, within the sliding window, the optimizer faces a prediction horizon that has shifted slightly but not significantly, and the difference between the high-precision Kriging surrogate model used in the optimization and the SPM model used by the simulator leads to only small control errors. The small discrepancies in control decisions lead to minor deviations in the SOC profile, but overall, the battery's state remains largely consistent.

However, after hour 10, the decisions made by the rolling optimization gradually begin to diverge more significantly from the initial planning. In some cases, entirely different decisions are made. For example, while the initial plan chooses to discharge the battery and not swap it at hour 11, the rolling optimization opts to charge the battery and perform the swap. This divergence leads to a noticeable shift in the SOC evolution. As the optimization moves further into the sliding window, the decisions regarding charging and swapping become increasingly different, resulting in a distinctly different trajectory for the SOC. This divergence occurs because, as the sliding window progresses further, the optimization problem becomes increasingly different from the original one. The new optimization problems must account for conditions outside the planned horizon, and the cumulative error from the Kriging surrogate model also starts to have a noticeable impact. This highlights how small changes in the optimization decisions can lead to significantly different battery behavior and SOC evolution over time.

In summary, BSS-MPC effectively responds to fluctuations in electricity prices and demand within the prediction horizon when solving the optimization problem. It accumulates a certain amount of electricity during low-price and low-demand periods, minimizing purchases during peak periods to save costs. Additionally, the system identifies opportunities for arbitrage by selling electricity during certain periods. Meanwhile, the rolling optimization process continuously updates the decision-making window, allowing for the adaptation of decisions as the considered prediction horizon changes. This approach also mitigates the cumulative errors at the end of the horizon caused by the Kriging surrogate model in one-time decision-making.

\subsection{A Half-Year Comparison of Three Methods}\label{sec:comparison}

To further compare the performance and characteristics of the three strategies, we evaluated their behavior over a 180-day simulation period under identical initial conditions.

\subsubsection{Profit and Capacity Fade Comparison}
\begin{figure*}[!h]
    \centering
    \includegraphics[width=1.0\linewidth]{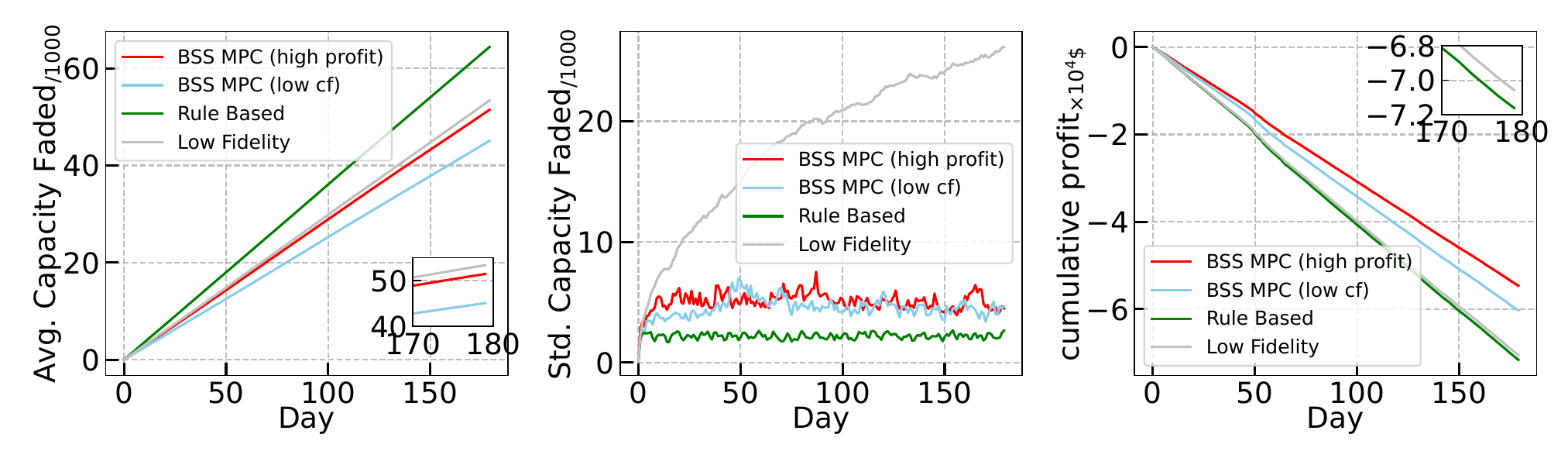}
    \caption{180-Day Cumulative Profit and Capacity Fade Comparison Among Strategies.}
    \label{fig:profit_and_cf_comparison}
\end{figure*}
The cumulative profit, average capacity degradation, and variance of \( c_f \) for the three strategies are summarized in ~\Cref{fig:profit_and_cf_comparison}.

The rule-based strategy exhibits the weakest overall performance. While selecting the battery with the largest $c_f$ yields the most uniform degradation, it completely ignores electricity prices and arbitrage opportunities. Instead, it blindly charges incoming batteries to full capacity within one hour, regardless of market conditions. This leads to continuous energy purchases without revenue and accelerates degradation due to aggressive charging. As a result, it incurs only costs—yielding the lowest profit—and suffers the highest average $c_f$.

Compared with the rule-based strategy, the low-fidelity model achieves a similar total profit while reducing the average degradation cost \( c_f \) by approximately \( 17\% \). This improvement mainly stems from grid-aware scheduling and the implicit power penalty, which captures the positive correlation between battery aging and charging/discharging power. By explicitly regularizing the power, the model encourages slower charging and discharging, thereby mitigating degradation.  

In addition, the low-fidelity model employs a simplified battery energy evolution model, which allows it to effectively account for variations in electricity prices and swapping demands. As a result, the model can reduce the rate of battery degradation without incurring significant additional charging or discharging costs.  

However, the variance of the degradation cost, is noticeably larger under the low-fidelity model, indicating uneven degradation across the battery fleet. This unevenness originates from the simplified dynamic model’s inability to account for the spatial distribution of capacity fade. Although the model captures high-level trade-offs between arbitrage and battery aging, it lacks the fine-grained control necessary to ensure uniform degradation among all batteries.

The \textbf{BSS-MPC} strategy demonstrates the highest flexibility and overall performance, evaluated under two configurations:

\begin{itemize}
    \item Under the \textit{high-profit} setting ($w_1 = 10^2$, $w_2 = 10^3$), it achieves the highest total profit (24\% improvement over rule-based), while also reducing average $c_f$ by 16\%.
    \item The \textit{low-$c_f$} setting ($w_1 = 10^3$, $w_2 = 5 \times 10^2$) yields the lowest average capacity fade (30\% improvement), with profit nearly matching that of the low-fidelity model (14\% improvement over rule-based).
\end{itemize}

Both versions demonstrate significant advantages in terms of profit and battery degradation compared with the rule-based and low-fidelity models. By employing a Kriging surrogate to approximate the high-fidelity battery model, the evolution of battery energy and capacity is accurately captured. This enables the BSS-MPC to precisely exploit arbitrage opportunities arising from variations in electricity prices and swapping demands. Meanwhile, the rate of battery aging is explicitly controlled, allowing the system to dynamically adjust parameters as needed to effectively balance short-term profit and long-term battery degradation.

In addition, by explicitly penalizing degradation distribution via ~\Cref{eq:distribution penalty}, BSS-MPC results in a $c_f$ variance approximately twice that of the rule-based strategy—still indicating reasonably balanced battery usage. Unlike the low-fidelity model, which passively balances profit and degradation through fixed model terms and exhibits significantly higher variance, BSS-MPC enables tunable control over the trade-off between economic gain and battery longevity. This adaptability allows BSS-MPC to either prioritize revenue or mitigate degradation, depending on application needs—achieving both lower average $c_f$ and more balanced degradation across the battery fleet.

\begin{figure}[!h]
    \centering
    \begin{subfigure}[b]{0.4\linewidth}
        \centering
        \includegraphics[width=\linewidth]{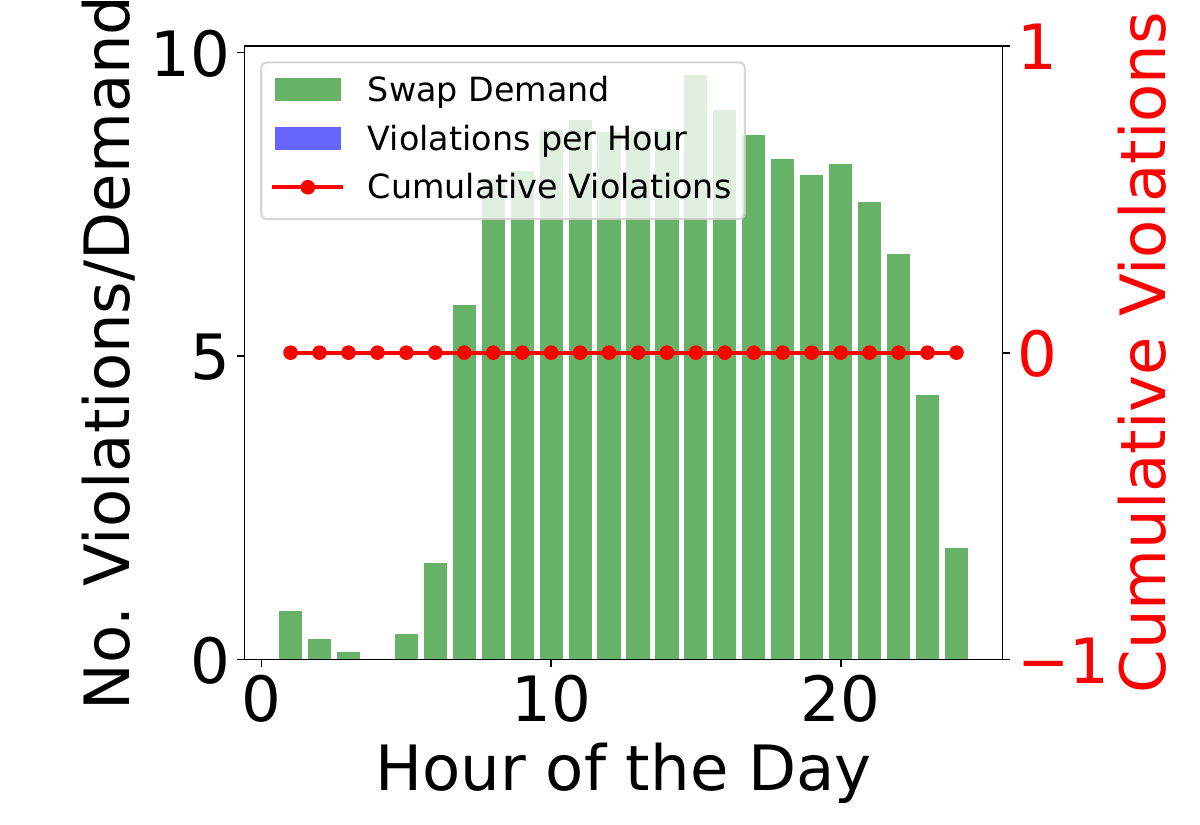}
        \caption{Rule Based.}
        \label{fig:violation_rule}
    \end{subfigure}
    \begin{subfigure}[b]{0.4\linewidth}
        \centering
        \includegraphics[width=\linewidth]{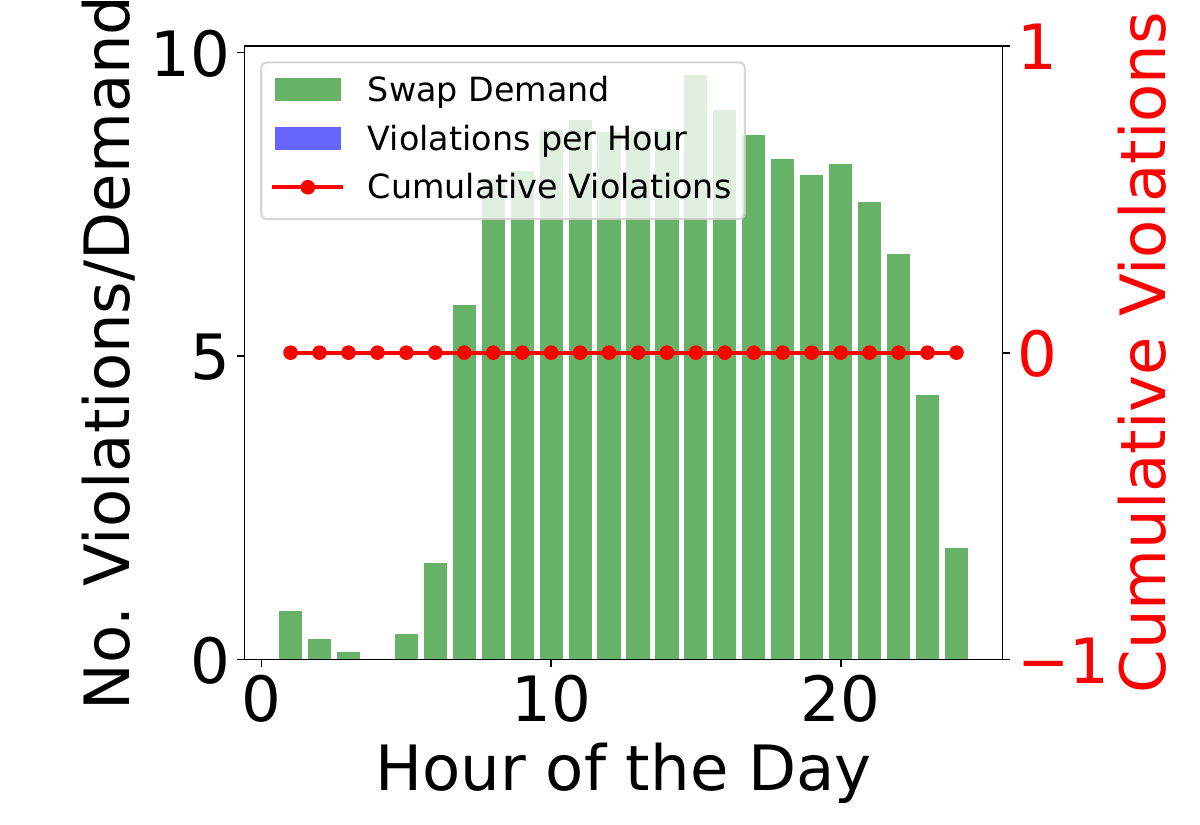}
        \caption{BSS-MPC (high-profit mode).}
        \label{fig:violation_high_profit}
    \end{subfigure}

    \begin{subfigure}[b]{0.4\linewidth}
        \centering
        \includegraphics[width=\linewidth]{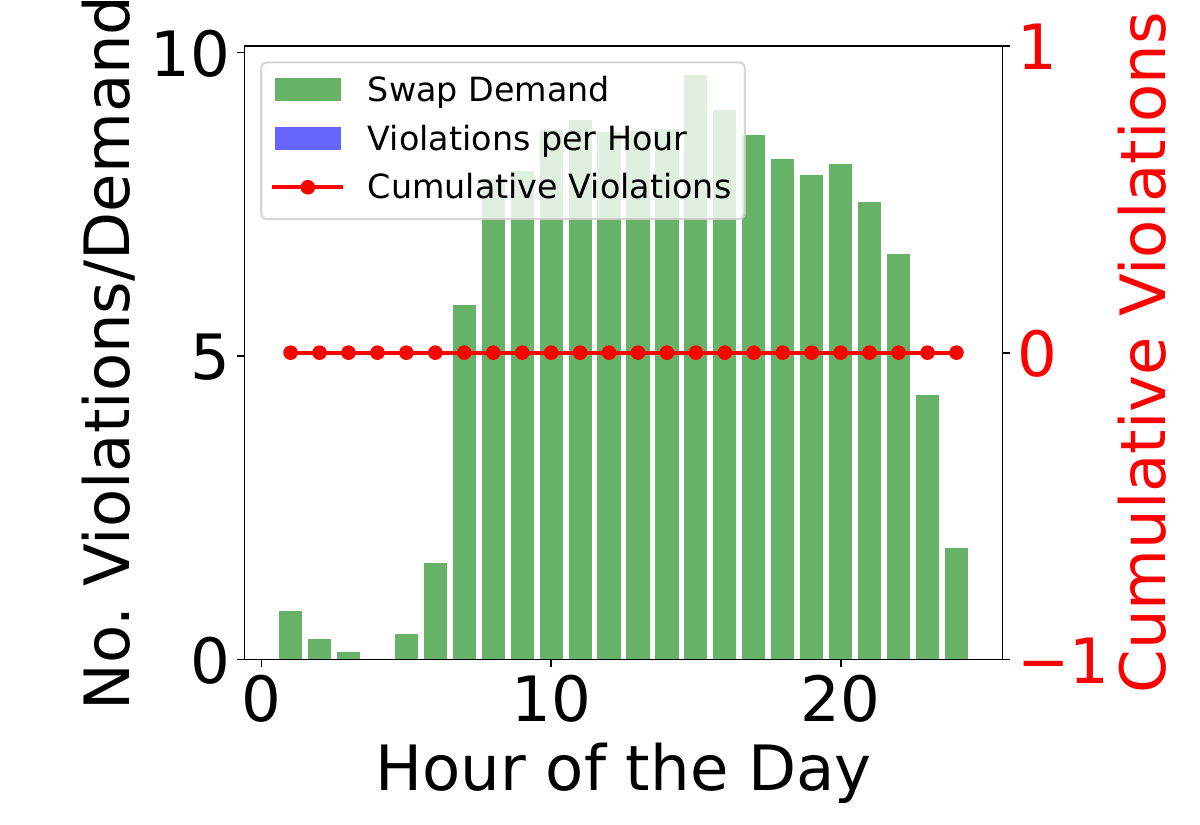}
        \caption{BSS-MPC (low-cf mode).}
        \label{fig:violation_low_cf}
    \end{subfigure}
    \begin{subfigure}[b]{0.4\linewidth}
        \centering
        \includegraphics[width=\linewidth]{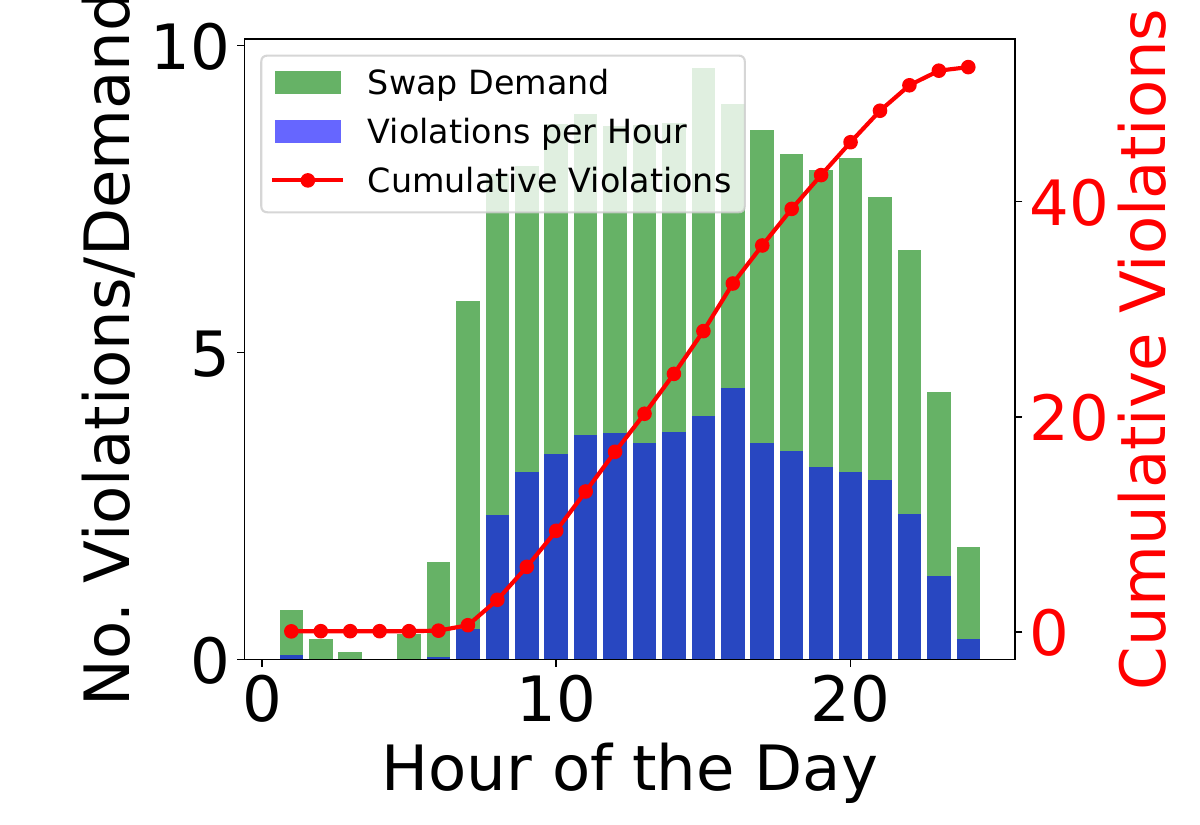}
        \caption{Low Fidelity Model.}
        \label{fig:violation_low_fidelity}
    \end{subfigure}

    \caption{Daily SOC-constraint violations across control strategies.}
    \label{fig:violation}
\end{figure}

\subsubsection{SOC Constraint Satisfaction}

All strategies incorporate a small tolerance $\varepsilon$ to accommodate model approximation and discretization errors. We recorded the average state-of-charge (SOC) constraint satisfaction over the 180-day period, measured by the number of batteries selected for swapping despite violating the SOC threshold (see ~\Cref{fig:violation}), where the threshold is set to \( \tilde{\mathrm{SOC}} = 0.7 \).

The low-fidelity model exhibits the poorest compliance. Even with a loose tolerance of $\varepsilon = 0.1$—two orders of magnitude larger than the others—it still results in frequent violations, with about 40\% of selected batteries falling below the threshold. Its oversimplified charge/discharge dynamics lead to SOC drifting both above 0.75 and below 0.7 (e.g., down to 0.65), making it difficult to eliminate violations solely by tuning $\varepsilon$.

In contrast, BSS-MPC uses a Kriging surrogate to accurately fit the true one-hour state transition. Since it directly models transitions at the hourly resolution, it avoids discretization errors altogether and only requires the relative fitting error to be sufficiently small. Despite using a much tighter tolerance ($\varepsilon = 0.001$), it satisfies all SOC constraints. The final SOC values are tightly clustered near the threshold, comparable to those produced by the high-fidelity model used in the dummy strategy. This highlights the high modeling accuracy of the Kriging-based dynamics, matching the reliability of rule-based control while enabling full optimization flexibility.

\subsubsection{Analysis of Cost Breakdown}
\begin{figure}[!h]
    \centering
    \includegraphics[width=0.8\linewidth]{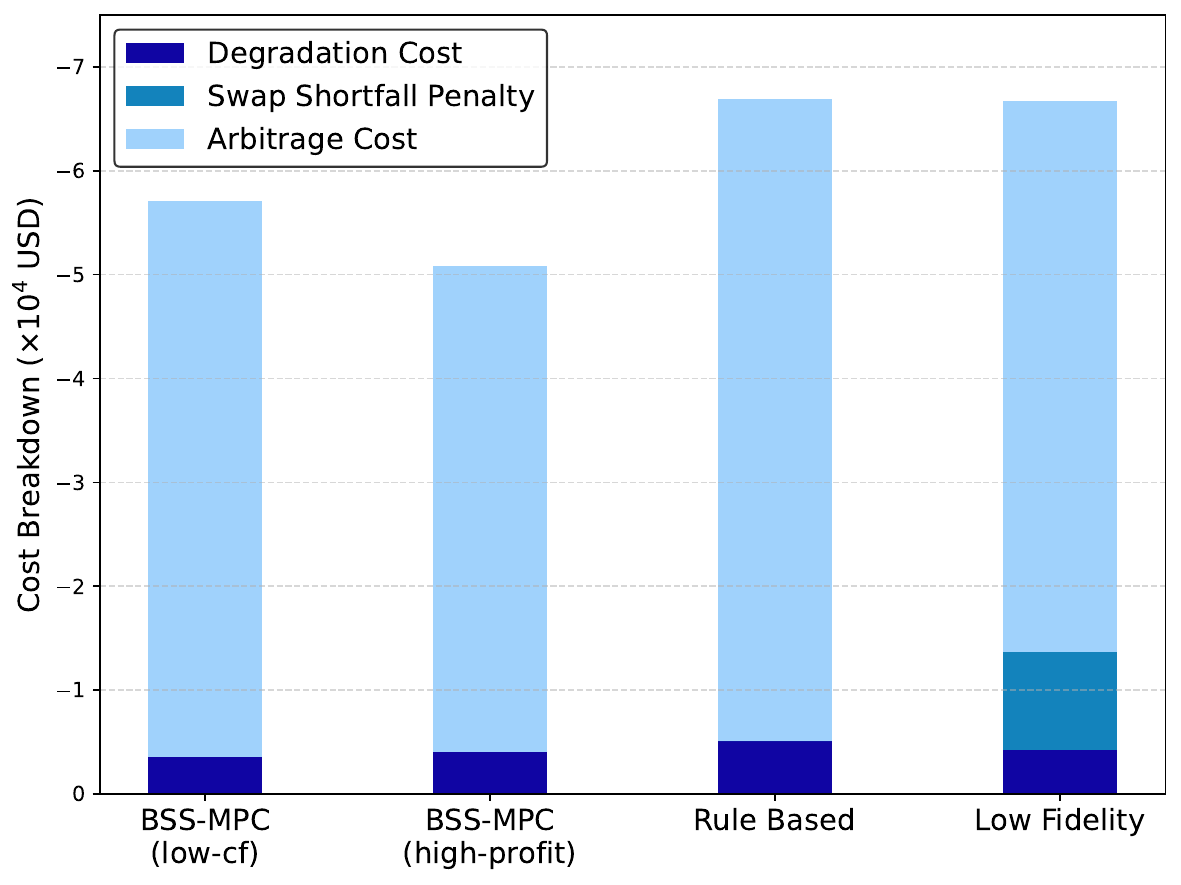}
    \caption{Cost breakdown of different control strategies}
    \label{fig:stacked_profit_by_method}
\end{figure}

Figure~\Cref{fig:stacked_profit_by_method} shows the cost breakdown of the three strategies. The total cost consists of three parts:  
(i) electricity cost 
(ii) depreciation cost due to battery aging, and  
(iii) penalties incurred when selected batteries do not meet the required SOC threshold.

Consistent with the observations in the previous sections, the high-profit version of BSS-MPC exhibits the greatest advantage in terms of charging and discharging costs. Although its battery degradation cost is slightly higher than that of the low-\( c_f \) version, the charging and discharging costs dominate, resulting in the lowest total cost. The low-\( c_f \) version, on the other hand, achieves the lowest aging depreciation cost while still maintaining a significant advantage in charging and discharging costs compared with the rule-based and low-fidelity models.  

Furthermore, both the BSS-MPC and the rule-based strategies strictly satisfy the SOC constraints, resulting in zero swapping penalties. In contrast, penalties account for a substantial portion of the total cost for the low-fidelity model. Although its charging and discharging costs are approximately \( 14\% \) lower than those of the rule-based strategy, frequent SOC violations incur significant penalties, which offset these gains and eliminate any cost advantage.  

In practical battery swapping station operations, the penalties associated with selecting batteries below the required SOC are likely to be far more severe than the simplified assumptions used in our experiments. Such violations can have long-term implications for the station’s reputation and operational sustainability, which are beyond the scope of this work to quantify.

\begin{table*}[!h]
\centering
\caption{Comparison of three strategies across four key metrics.}
\label{tab:strategy_comparison}
\begin{tabular}{lcccc}
\toprule
\textbf{Strategy} & \textbf{Normalized Loss (\%)} & \textbf{Avg $c_f$ (\%)} & 
\textbf{$c_f$ Variance (\%)} & \textbf{SOC Satisfaction (\%)} \\
\midrule
Rule-Based & 100 & 100 & \textbf{100} & \textbf{100} \\
Low-Fidelity & 98.6 & 82.79 & 1160.7 & 60 \\
BSS-MPC (high-profit) & \textbf{76.04} & 79.92 & 176.14 & \textbf{100} \\
BSS-MPC (low-cf) & 85.29 & \textbf{70.05} & 212.63 & \textbf{100} \\
\bottomrule
\end{tabular}
\end{table*}
\subsubsection{Overall Comparison}
To summarize, ~\Cref{tab:strategy_comparison} compares the three strategies across four key metrics: cumulative profit, average $c_f$, $c_f$ variance, and SOC compliance. All values are normalized to the rule-based strategy (set as 100\%) for intuitive comparison, with the best performance in each metric highlighted in bold. Since all absolute profit values are negative (indicating losses), we present \textbf{Normalized Loss} instead of profit, defined relative to the rule-based strategy. Here, a lower Normalized Loss percentage indicates better economic performance (i.e., smaller cumulative losses).

As shown in ~\Cref{tab:strategy_comparison}, \textbf{BSS-MPC} consistently excels across all four evaluation metrics, demonstrating its balanced and adaptable design.

\begin{itemize}
    \item \textbf{Economic performance}: Under the high-profit configuration, BSS-MPC results in the \textbf{least financial loss (76.04\% normalized profit)}, outperforming both the low-fidelity model and rule-based baseline. Even in the low-$c_f$ setting, it maintains relatively low loss (85.29\%), striking an effective trade-off between revenue and degradation.
    
    \item \textbf{Degradation control}: BSS-MPC significantly reduces \textbf{average capacity fade}, with the low-$c_f$ configuration achieving the \textbf{lowest $c_f$ (70.05\%)}. This indicates a \textbf{notable improvement in battery longevity}. Furthermore, both BSS-MPC variants maintain much lower \textbf{$c_f$ variance} compared to the low-fidelity model (e.g., 176.14\% vs. 1160.7\%), implying \textbf{more uniform battery usage} across the fleet and fewer aging outliers.
    
    \item \textbf{SOC constraint satisfaction}: Thanks to Kriging-based modeling that directly fits one-hour transitions without discretization error, BSS-MPC enforces \textbf{strict SOC compliance (100\%)} even under tight tolerance $\varepsilon=0.001$. This accuracy matches the rule-based strategy while retaining full optimization flexibility---a capability the low-fidelity model lacks.
\end{itemize}

In contrast, the \textbf{low-fidelity model}, while incurring slightly less loss than the rule-based baseline, suffers from \textbf{high $c_f$ variance} and \textbf{low SOC compliance}, indicating \textbf{inefficient battery utilization and potential stability concerns} in operation. The \textbf{rule-based strategy}, although fully compliant and uniform in battery usage by design, is \textbf{overly conservative}, leading to the \textbf{highest financial loss}.

In summary, BSS-MPC achieves a rare combination of economic viability, degradation mitigation, and robust constraint enforcement, proving itself the most capable and versatile strategy among the three. All the numerical experiments were conducted using warm-start. Under this setup, each simulation iteration of BSS-MPC takes approximately 20–40 seconds to complete, demonstrating the its capability for near real-time operation.

\section{Conclusion}
This paper presents \emph{BSS-MPC}, the first real-time MPC framework for BSS that explicitly embeds a physics-informed capacity fade model. By capturing both fast-timescale state-of-charge evolution and slow-timescale degradation dynamics, the framework enables operators to systematically balance short-term energy arbitrage profit with long-term battery lifetime, outperforming conventional rule-based and low-fidelity scheduling approaches. \emph{BSS-MPC} offers a practical and scalable solution for sustainable EV infrastructure. Future work will focus on further reducing computational complexity by developing faster surrogate models and exploring stochastic MPC to handle uncertainties in electricity prices and swapping demand. 

\section*{Acknowledgments}

\bibliography{main}
\end{document}